\documentclass[10pt]{article}

\usepackage{amsthm}
\usepackage{amsmath,amssymb}
\usepackage[dvips]{graphicx}
\usepackage{float}
\usepackage{wrapfig}
\usepackage{subcaption}
\usepackage[normalem]{ulem}

\newtheorem{theorem}{Theorem}[section]

\theoremstyle{definition}
\newtheorem{definition}[theorem]{Definition}

\newtheorem{example}[theorem]{Example}
\newtheorem{corollary}[theorem]{Corollary}

\theoremstyle{remark}
\newtheorem{remark}[theorem]{Remark}

\title{On link diagrams that are minimal with respect to Reidemeister moves I and II}
\author{Kishin Sasaki}

\begin{document}
\maketitle
\begin{abstract}
In this paper, a link diagram is said to be minimal if no Reidemeister move I or II can be applied to it to reduce the number of crossings. We show that for an arbitrary diagram $D$ of a link without a trivial split component, a minimal diagram obtained by applying Reidemeister moves I and II to $D$ is unique. 
The proof also shows that the number of crossings of such a minimal diagram is unique for any diagram of any link.
We show that for a link without a trivial split component, an arbitrary Reidemeister move III either does not change the associated minimal diagram or can be reduced to a special type of a move up to Reidemeister moves I and II.
\end{abstract}

\section{Introduction}
It is well known that every pair of diagrams of a given link can be transformed to each other by applying finitely many Reidemeister moves. Furthermore, any two diagrams of links that are transformed to each other by finitely many Reidemeister moves represent equivalent links (\cite{4}). There have been a lot of studies about the Reidemeister moves (see, for example, \cite{1,2,3,6,7}).

In this paper, we consider smooth unoriented link diagrams in $\mathbb{R}^{2}$ or $S^{2}$. In Section~\ref{preliminaries}, we prepare some terminologies and notions necessary for later sections. Reidemeister moves I and II change the number of crossings, while Reidemeister move III does not change the number of crossings. From this viewpoint, in Section~\ref{sec:minimum}, we say that a link diagram is minimal if no Reidemeister move I or II can be applied to it to reduce the number of crossings, and show that for an arbitrary diagram $D$ of a link without a trivial split component, a minimal diagram obtained by applying Reidemeister moves I and II to $D$ is unique (Theorem~\ref{theo:minimality}). Furthermore, the proof also shows that the number of crossings of such a minimal diagram is uniquely determined, for a diagram of an arbitrary link possibly with a trivial split component (Corollary~\ref{cor:minimal_crossing_number}). The idea of this minimality has appeared in [\cite{5}, Theorem 2.2], where the uniqueness of the minimality of knot projections without crossings' information has been studied by Mikhail Khovanov.

By studying the Reidemeister move III from the viewpoint of minimal diagram change, in Section~\ref{sec:relation},
we show that for a link without a trivial split component, an arbitrary Reidemeister move III either does not change the associated minimal diagram or can be reduced to a special type of a move up to Reidemeister moves I and II (Theorem~\ref{theo:RIII_reduced}). As a corollary, we will see that, for every RI-II equivalence class (see Definition~\ref{def:RI-II equivalence}) of a link without a trivial split component, the set of RI-II equivalence classes $(-)$-adjacent (see Definition \ref{def:adjacent}) to the original RI-II equivalence class is the set of RI-II equivalence classes obtained by applying a Reidemeister move III or III* to the minimal diagram in the original RI-II equivalence class (Corollary~\ref{cor:equivalence_relation}). The result enhances the utility of the minimal diagrams.

\section{Preliminaries}
\label{preliminaries}
In this section, we give definitions and a remark which will be used in Sections~\ref{sec:minimum}~and~\ref{sec:relation}.

\begin{figure}[H]
\centering
\begin{subfigure}[t]{0.4\textwidth}
\centering
\includegraphics[width=\textwidth]{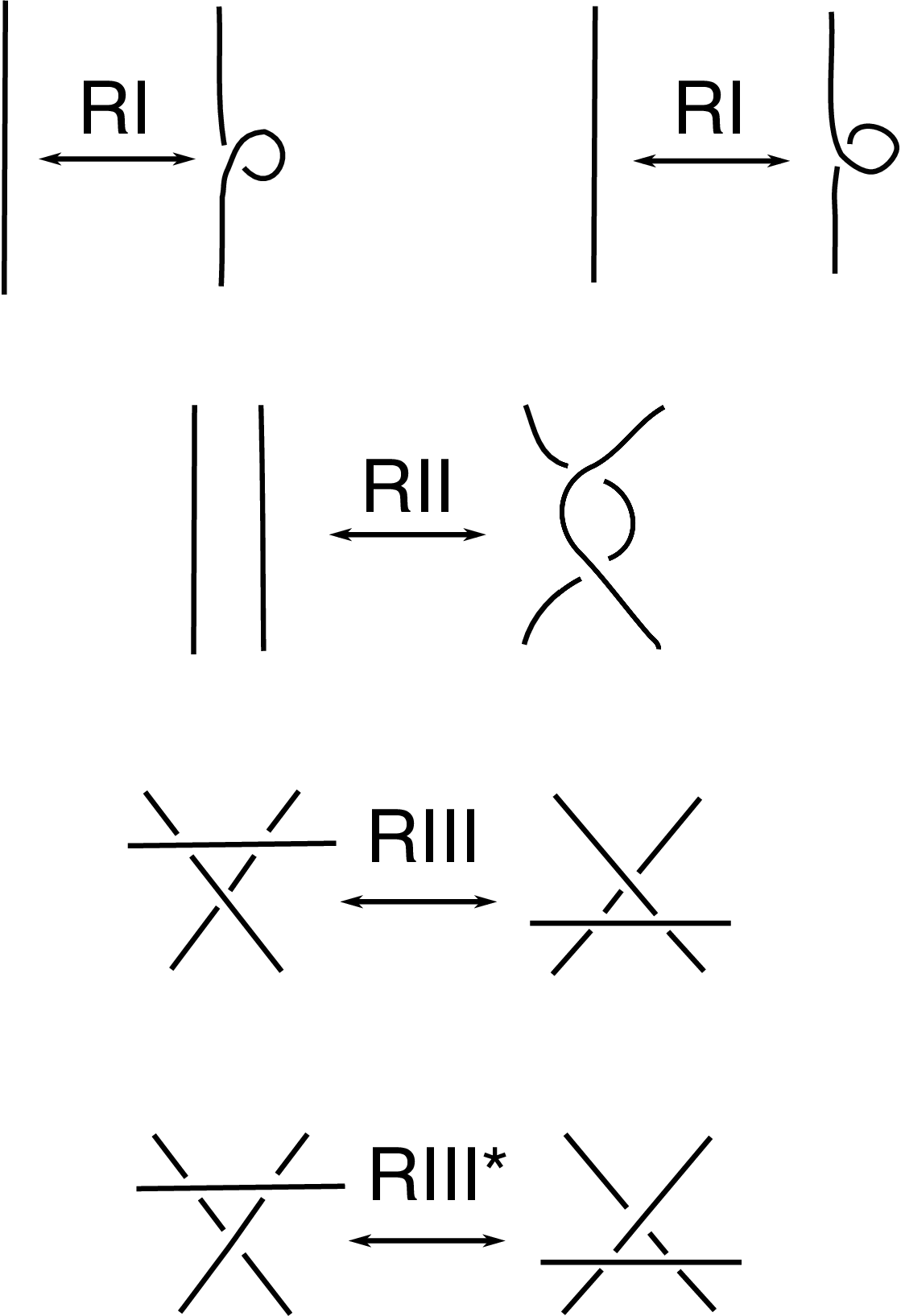}
\caption{Reidemeiter moves I, II, III}
\label{fig:Reidemeister}
\end{subfigure}
\hfill
\begin{subfigure}[t]{0.5\textwidth}
\centering
\includegraphics[width=\textwidth]{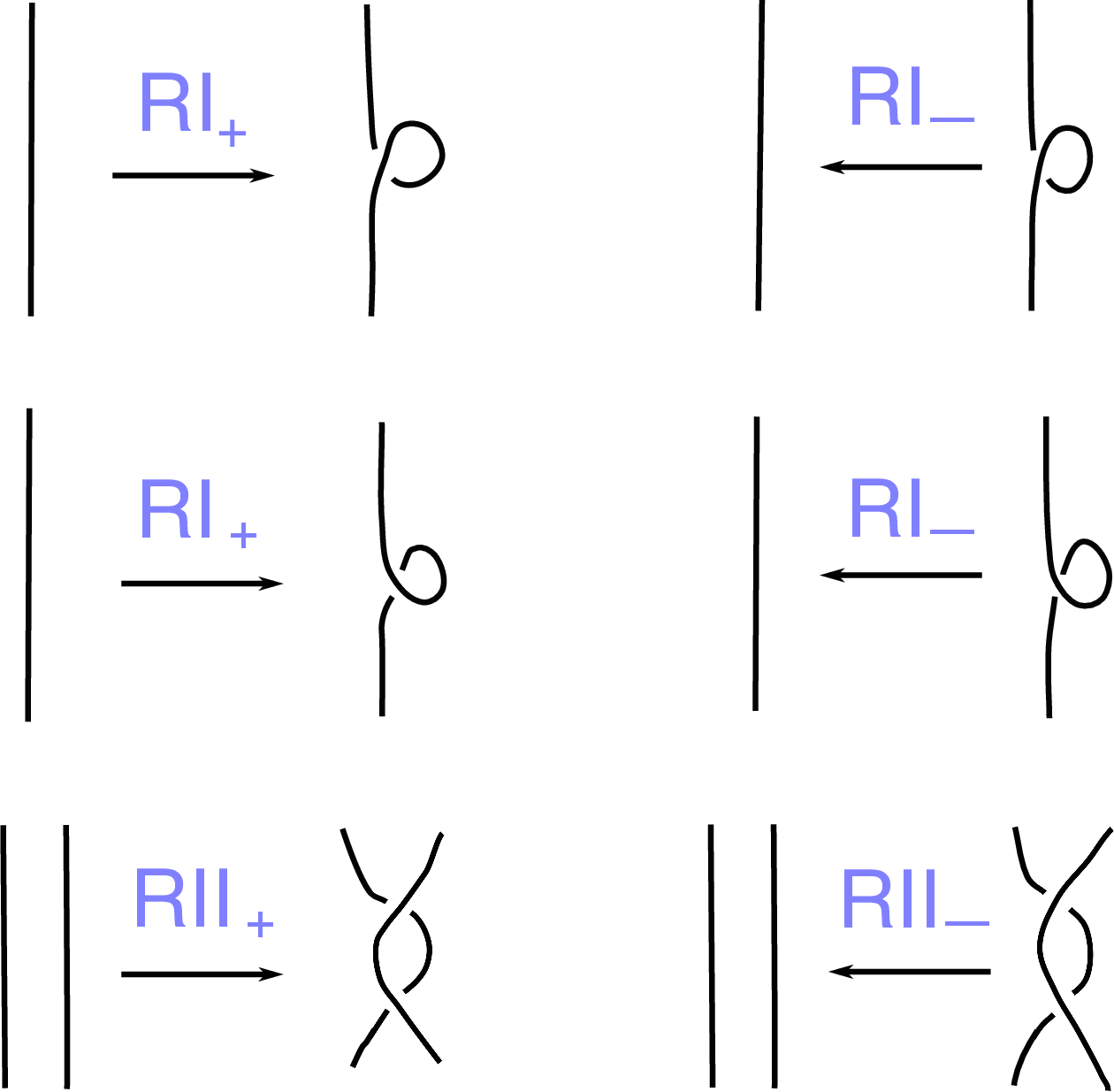}
\caption{The moves $RI_{+}$, $RII_{+}$, $RI_{-}$, $RII_{-}$}
\label{fig:RI+-,RII+-}
\end{subfigure}
\caption{}
\end{figure}

\begin{definition}
\emph{Reidemeister moves} are defined as the local moves of link diagrams as depicted in Figure \ref{fig:Reidemeister}. The moves RI, RII, RIII and RIII* depicted in Figure \ref{fig:Reidemeister} are called Reidemeister moves I, II, III and III*, respectively. The moves RIII and RIII* can be distinguished by using the orientation of $\mathbb{R}^{2}$ (or $S^{2}$). 
\end{definition}

The following theorem is well known.

\begin{theorem}[Reidemeister \cite{4}, 1927]
Every pair of diagrams of a link may be transformed to each other by applying finitely many Reidemeister moves. Furthermore, any two diagrams of links which are transformed to each other by applying finitely many Reidemeister moves represent equivalent links.    
\end{theorem}

\begin{definition}
Moves $\mathrm{RI}_+$ and $\mathrm{RII}_+$ are defined to be Reidemeister moves I and II which increase the number of crossings, respectively. Moves $\mathrm{RI}_-$ and $\mathrm{RII}_-$ are defined to be Reidemeister moves I and II which decrease the number of crossings, respectively. See Figure \ref{fig:RI+-,RII+-}. 
\end{definition}

\begin{remark}
\label{rem:shape}
Reidemeister moves I, II, III and III* are related to a monogon, the digon, and the two triangles, respectively, which appear in the local disks on the right hand sides of Figure \ref{fig:Reidemeister} (the triangles also appear in the local disks on the left hand sides of Reidemeister moves III, III*). 
\end{remark}

\section{Minimal link diagrams with respect to Reidemeister moves I and II}
\label{sec:minimum}
In this section, for every diagram of a link without a trivial split component, we prove that all the diagrams obtained by applying finitely many moves $\mathrm{RI}_-$ and $\mathrm{RII}_-$ until they cannot be applied are equivalent (or planer isotopic). We also show that for an arbitrary diagram of a link possibly with a trivial split component, the number of crossings of such a minimal diagram is uniquely determined. 

\begin{definition}
\label{def:minimal}
A link diagram is said to be \emph{minimal} in this paper if no Reidemeister move $\mathrm{RI}_-$ or $\mathrm{RII}_-$ can be applied to it.
\end{definition}

\begin{theorem}\label{theo:minimality}
For an arbitrary diagram $D$ of a link without a split unknot component, all the minimal diagrams which are obtained by applying finitely many moves $\mathrm{RI}_-$ and $\mathrm{RII}_-$ to $D$ are equivalent. 
\end{theorem}

\begin{remark}

[\cite{5},Theorem 2.2] gives a result similar to Theorem \ref{theo:minimality}. Mikhail Khovanov has studied the uniqueness of  the minimality of knot projections without crossings' information there. 
\end{remark}

\begin{proof}
Let $L$ be a link without a split unknot component. From a given diagram $D$ of this link $L$, we clearly obtain a minimal diagram
$D_{1}$ by applying finitely many moves $\mathrm{RI}_-$ and $\mathrm{RII}_-$. Note that this minimal diagram $D_{1}$ may depend on the sequence of the moves $\mathrm{RI}_-$ and $\mathrm{RII}_-$ applied to the diagram $D$. The minimal diagram $D_{1}$ clearly satisfies condition (\#) below.

\medskip

(\#)\quad All minimal diagrams obtained by applying finitely many moves $\mathrm{RI}_-$ and $\mathrm{RII}_-$ to that diagram are equivalent.

\medskip

Let us prove that all the diagrams which are obtained by applying finitely many moves $\mathrm{RI}$ and $\mathrm{RII}$ to the minimal diagram $D_{1}$ satisfy condition (\#) by checking that moves RI and RII maintain condition (\#).

First, we prove that moves RI maintain condition (\#).

Moves $\mathrm{RI}_-$ clearly maintain condition (\#) because of the definition of condition (\#).

Let us now prove the case of moves $\mathrm{RI}_+$. In other words, we prove that, as in Figure~\ref{fig:RI+}, the diagram on the right hand side satisfies condition (\#) if we assume that the diagram on the left hand side satisfies condition (\#).

\begin{figure}[H]
\centering
\includegraphics[width=.45\linewidth]{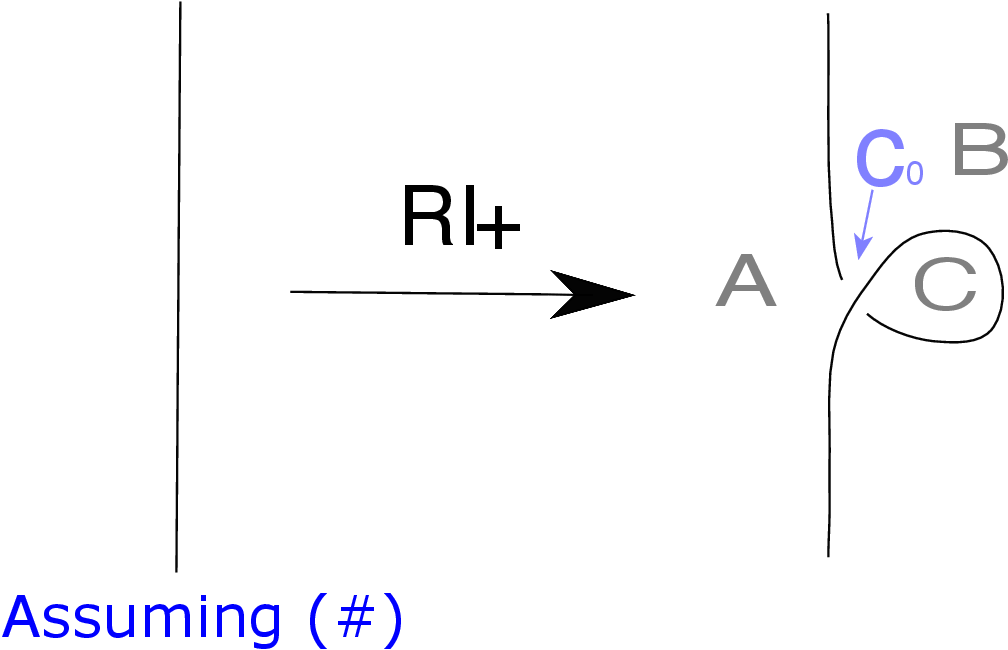}
\caption{Applying $\mathrm{RI}_+$}
\label{fig:RI+}
\end{figure}

We consider a process of applying finitely many moves $\mathrm{RI}_-$ and $\mathrm{RII}_-$ to the diagram on the right hand side in Figure~\ref{fig:RI+} until we get a minimal diagram. If this process does not eliminate crossing $c_0$ inside the local disk, then the resulting diagram inside this local disk remains the same after this process, which is a contradiction. So this process necessarily eliminates crossing $c_0$. If the crossing(s) eliminated by the first move in this process are outside the local disk, then the same move corresponding to the same position can be applied to the diagram on the left hand side in Figure \ref{fig:RI+}, and the resulting diagram still satisfies condition (\#). This maintains the relation between the diagrams on the both hands sides in Figure \ref{fig:RI+} in terms of condition (\#). Let us now consider the first move $m_1$ in the process that eliminates the crossing $c_0$. Note that whenever the crossing $c_0$ is eliminated by the move $m_1$, region $A$, $B$ or $C$ is a monogon or the corresponding digon (see Remark~\ref{rem:shape}).

When $m_1$ is the move $\mathrm{RI}_-$ that uses region $A$ as a monogon, as in Figure \ref{fig:RI+1}, the diagram which is obtained by applying the move $m_1$ is equivalent to the diagram on the left hand side which satisfies condition (\#). So in this case, the diagram on the right hand side in Figure \ref{fig:RI+1} satisfies condition (\#).

\begin{figure}[H]
\centering
\includegraphics[width=.45\linewidth]{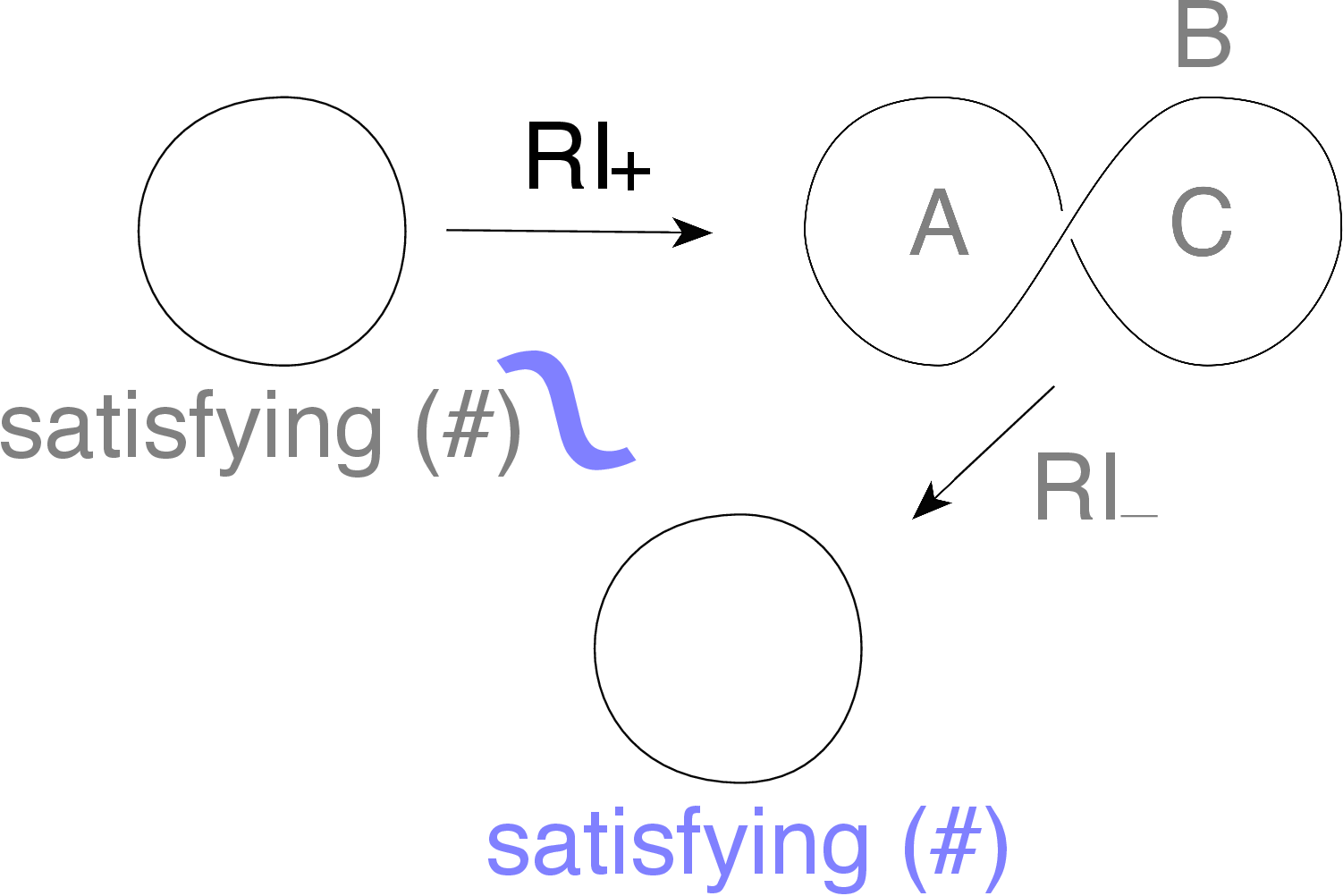}
\caption{The case of applying $\mathrm{RI}_-$ to region A}
\label{fig:RI+1}
\end{figure}

When $m_1$ is the move $\mathrm{RII}_-$ that uses region $A$ as the digon, as in Figure \ref{fig:RI+2}, the diagram which is obtained by applying the move $m_1$ is obtained by applying a move $\mathrm{RI}_-$ to the diagram on the left hand side, which satisfies condition (\#).
Also in this case, the diagram on the right hand side in Figure \ref{fig:RI+2} satisfies condition (\#).

\begin{figure}[H]
\centering
\includegraphics[width=.45\linewidth]{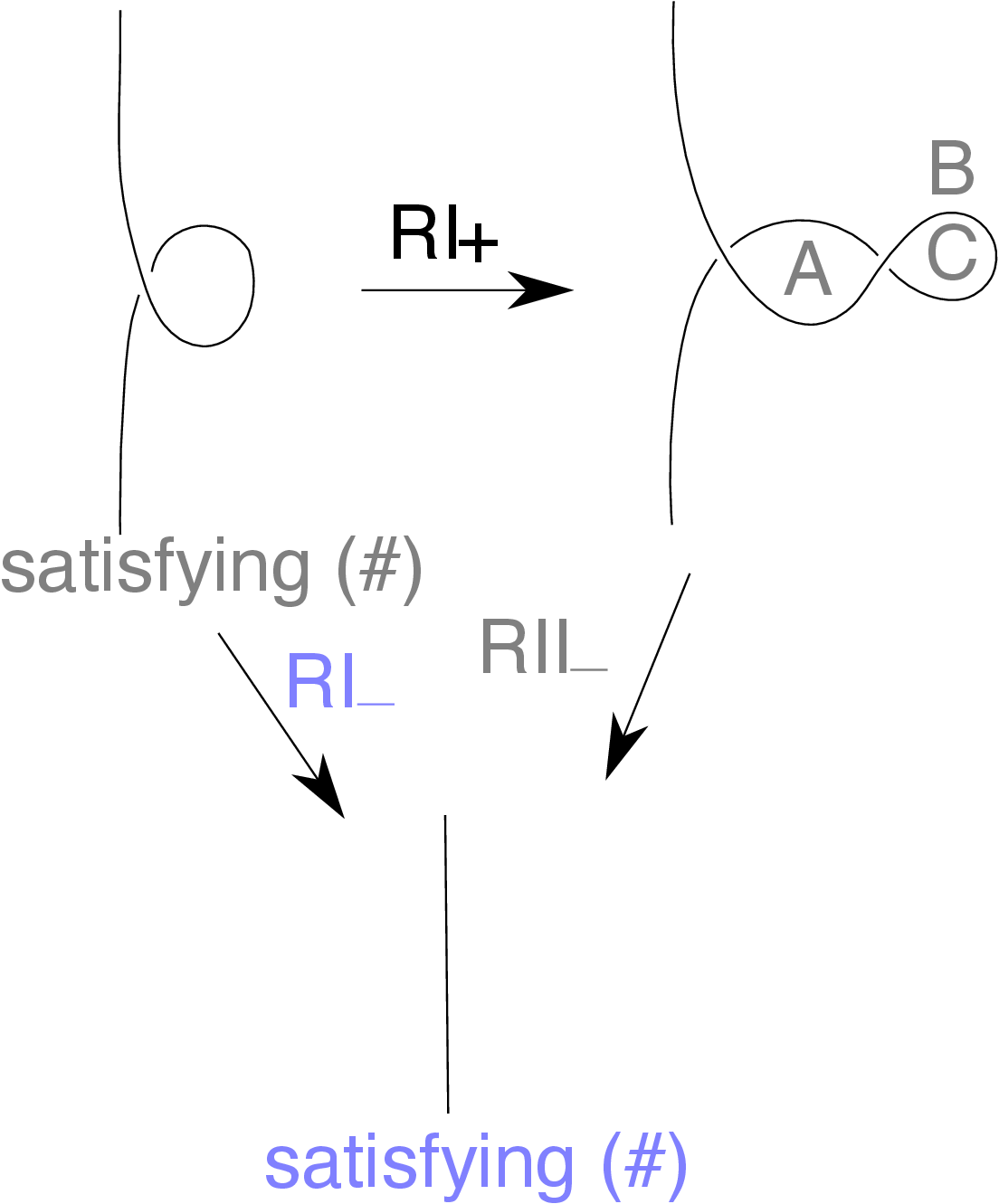}
\caption{The case of applying $\mathrm{RII}_-$ to region A}
\label{fig:RI+2}
\end{figure}

It is easy to see that region $B$ cannot be a monogon or the digon.

When $m_1$ is the move $\mathrm{RI}_-$ that uses region $C$ as a monogon, the diagram on the right hand side in Figure \ref{fig:RI+} clearly satisfies condition (\#), since the diagram on the right hand side which is obtained by applying the move $m_1$ is equivalent to the diagram on the left hand side which satisfies condition (\#).

Thus, in all the cases for moves $\mathrm{RI}_{+}$, we have checked that also the diagram on the right hand side in Figure \ref{fig:RI+} satisfies condition (\#), so moves $\mathrm{RI}_{+}$ maintain condition (\#).

Let us now prove that moves $\mathrm{RII}_{+}$ maintain condition (\#). Let us check all the cases for moves RII as we have checked all the cases for moves RI above. Moves $\mathrm{RII}_-$ clearly maintains condition (\#) because of the definition of condition (\#).

So let us now prove that moves $\mathrm{RII}_+$ maintain condition (\#).
See Figure~\ref{fig:RII+} which depicts a move $\mathrm{RII}_+$ whose regions adjacent to the two crossings inside the local disk on the right hand side are indicated by alphabets, and we assume that the diagram on the left hand side satisfies condition (\#). 

\begin{figure}[H]
\centering
\includegraphics[width=.45\linewidth]{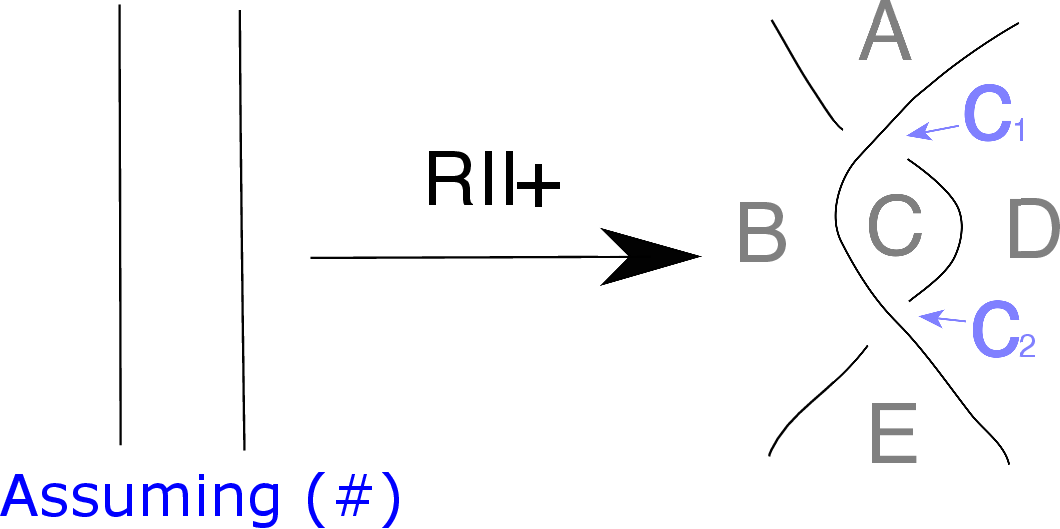}
\caption{Applying $\mathrm{RII}_+$}
\label{fig:RII+}
\end{figure}

We consider a process of applying finitely many moves $\mathrm{RI}_-$ and $\mathrm{RII}_-$ to the diagram on the right hand side in Figure~\ref{fig:RII+} until we get a minimal diagram. If this process does not eliminate a crossing inside the local disk (the upper side (or lower side) crossing is said to be $c_1$ (resp. $c_2$)), then the resulting diagram inside this local disk remains the same after this process, which is a contradiction. So this process necessarily eliminates crossing $c_1$ or $c_2$. If crossing(s) eliminated by the first move in this process are outside the local disk, then the same move corresponding to the same position can be applied to the diagram on the left hand side in Figure \ref{fig:RII+} and the resulting diagram still satisfies condition (\#). This maintains the relation between the diagrams on the left and the right hand sides in Figure \ref{fig:RII+} in terms of condition (\#). Let us now consider the first move $l_1$ in the process that eliminates crossing $c_1$ or $c_2$.

When $l_1$ is the move $\mathrm{RI}_-$ that uses the region $A$ as a monogon,
as in Figure \ref{fig:RII+1}, the diagram which is obtained by applying the move $l_1$ is transformed from the diagram on the left hand side by applying a move $\mathrm{RI}_{+}$, which maintains condition (\#) (we have proved this above). So in this case, also the diagram on the right hand side in Figure \ref{fig:RII+} satisfies condition (\#).

\begin{figure}[H]
\centering
\includegraphics[width=.45\linewidth]{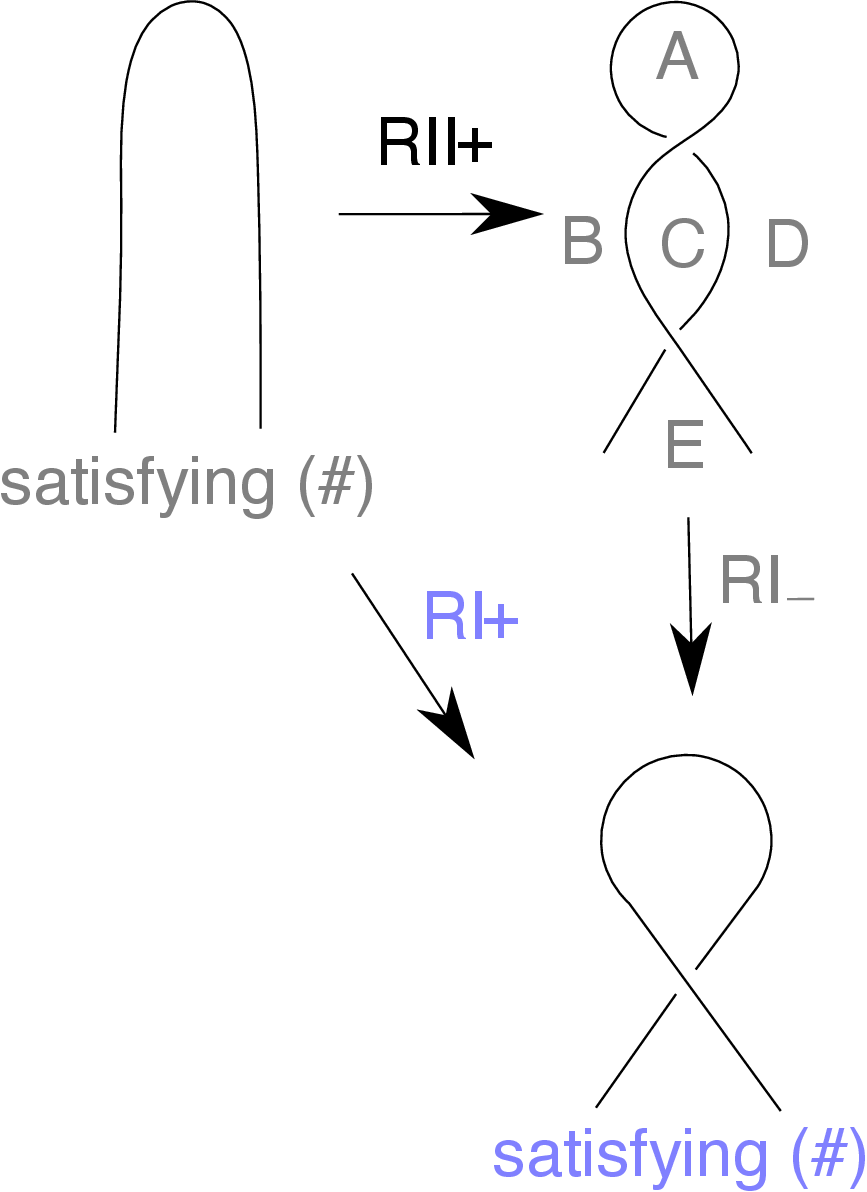}
\caption{Applying $\mathrm{RI}_-$ to region A}
\label{fig:RII+1}
\end{figure}

When $l_1$ is the move $\mathrm{RI}_-$ that uses the region $E$ as a monogon, by looking at Figure \ref{fig:RII+1} upside down, we can easily prove that also the diagram on the right hand side in Figure \ref{fig:RII+} satisfies condition (\#). 

When $l_1$ is the move $\mathrm{RII}_-$ that uses the region $A$ as the digon, as in Figure~\ref{fig:RII+2}, the diagram which is obtained by applying the move $l_1$ is equivalent to the diagram on the left hand side which satisfies condition (\#). So in this case, also the diagram on the right hand side in Figure \ref{fig:RII+} satisfies condition (\#).

\begin{figure}[H]
\centering
\includegraphics[width=.45\linewidth]{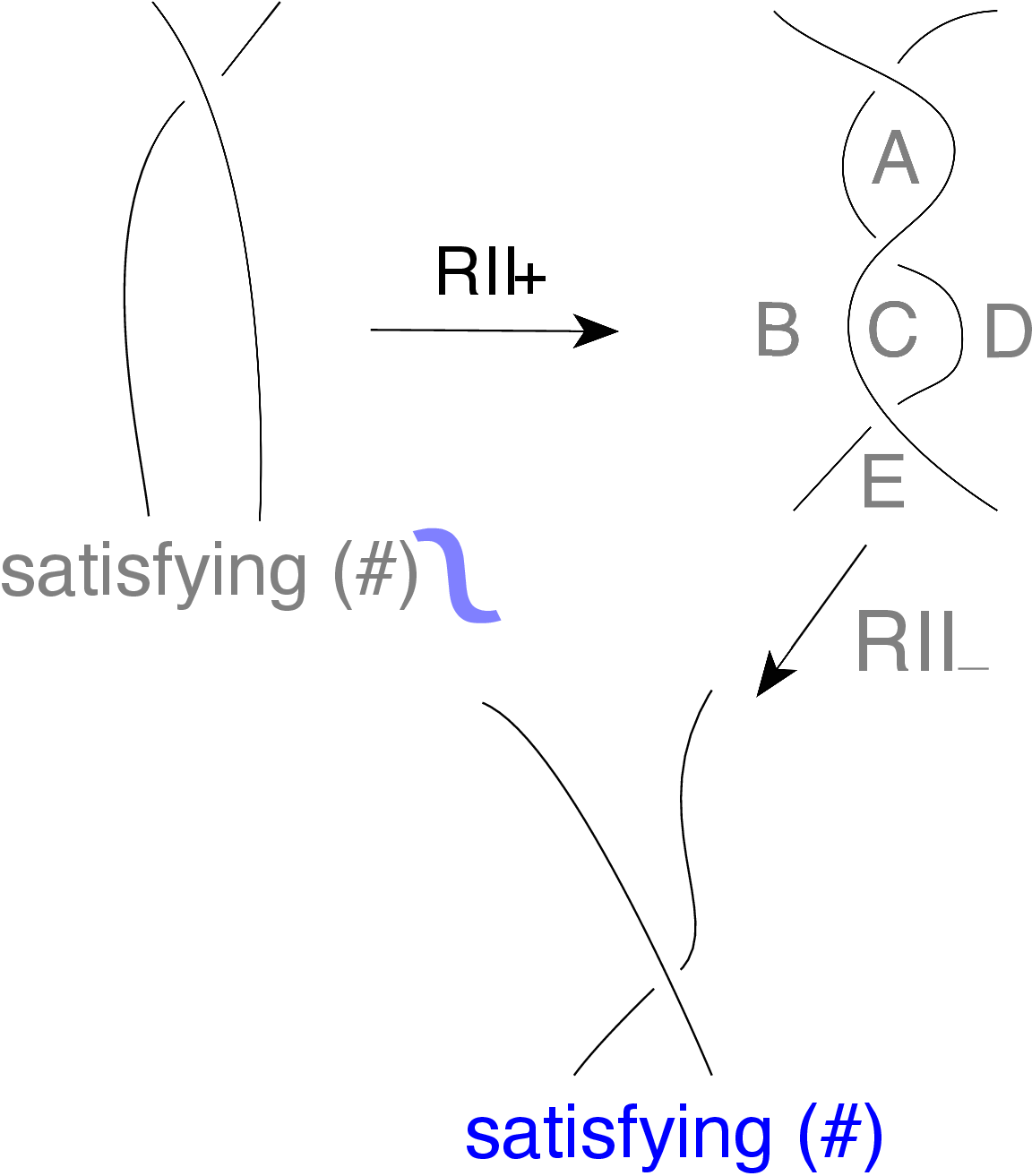}
\caption{Applying $\mathrm{RII}_-$ to region $A$}
\label{fig:RII+2}
\end{figure}

When $l_1$ is the move $\mathrm{RII}_-$ that uses the region $E$ as the digon, the proof is the same as the case of when $l_1$ is the move $\mathrm{RII}_-$ that uses the region $A$ as the digon.

When $l_1$ is the move $\mathrm{RII}_-$ that uses the region $C$ as the digon, also the diagram on the right hand side clearly satisfies condition (\#).

By the assumption of Theorem \ref{theo:minimality}, link $L$ does not have a split trivial component, so the case that $l_1$ is the move $\mathrm{RII}_-$ that uses region $B$ or $C$ as the digon cannot occur.

Thus, all the cases for moves $\mathrm{RII}_+$ have been checked. In all the cases, moves $\mathrm{RII}_+$ maintain condition (\#). We have proved that moves RI and RII maintain condition (\#) above, which completes the proof of Theorem~\ref{theo:minimality}
\end{proof}

\begin{example}

Figure~\ref{fig:counter_example} depicts a counterexample of the link-version Theorem \ref{theo:minimality}. We can get inequivalent minimal link diagrams by applying $\mathrm{RII}_{-}$ moves to the diagram on the left hand side in Figure \ref{fig:counter_example}.

\begin{figure}[H]
\centering
\includegraphics[width=.6\linewidth,keepaspectratio]{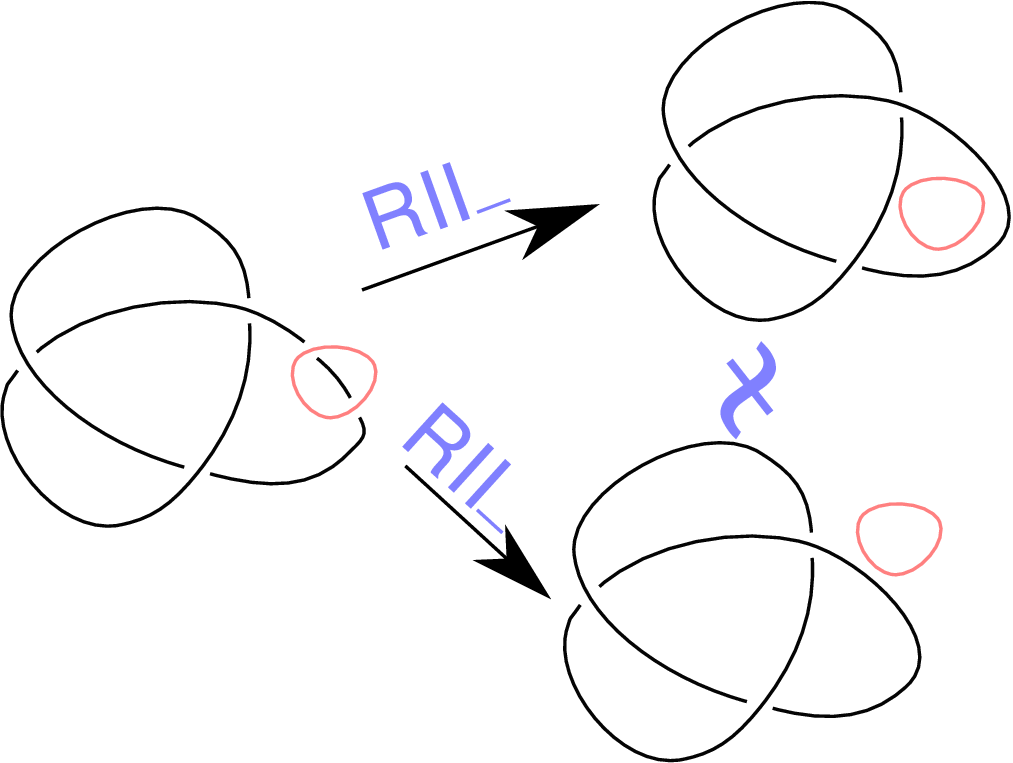}
\caption{A counter example for the link-version Theorem~\ref{theo:minimality}}
\label{fig:counter_example}
\end{figure}

\end{example}

\begin{corollary}
\label{cor:minimal_crossing_number}
For every diagram $D$ of a link possibly with trivial split components, we can get minimal diagrams with the same numbers of crossings by applying finitely many moves $\mathrm{RI}_-$ and $\mathrm{RII}_-$ to $D$.
 
\end{corollary}

\begin{proof}
By Theorem~\ref{theo:minimality}, for every diagram of a link without a trivial split component, the number of crossings of the minimal diagram is unique. For every diagram of a link with trivial split components, by the proof of Theorem~\ref{theo:minimality}, all the differences in the minimal diagrams are the locations of trivial split components which contain no crossing. So the numbers of crossings of the minimal diagrams are the same, which completes the proof.
\end{proof}

\begin{remark}
Figure~\ref{minimal_diagram} depicts an example of minimal diagrams of the trivial knot. It is easy to check that this diagram is minimal, since there are no monogons and no special digons for Reidemeister move II in this diagram. By considering connected sums of copies of this diagram and a minimal diagram, we can see that an arbitrary link has infinitely many minimal diagrams.

\begin{figure}[H]
\centering
\includegraphics[width=.3\linewidth]{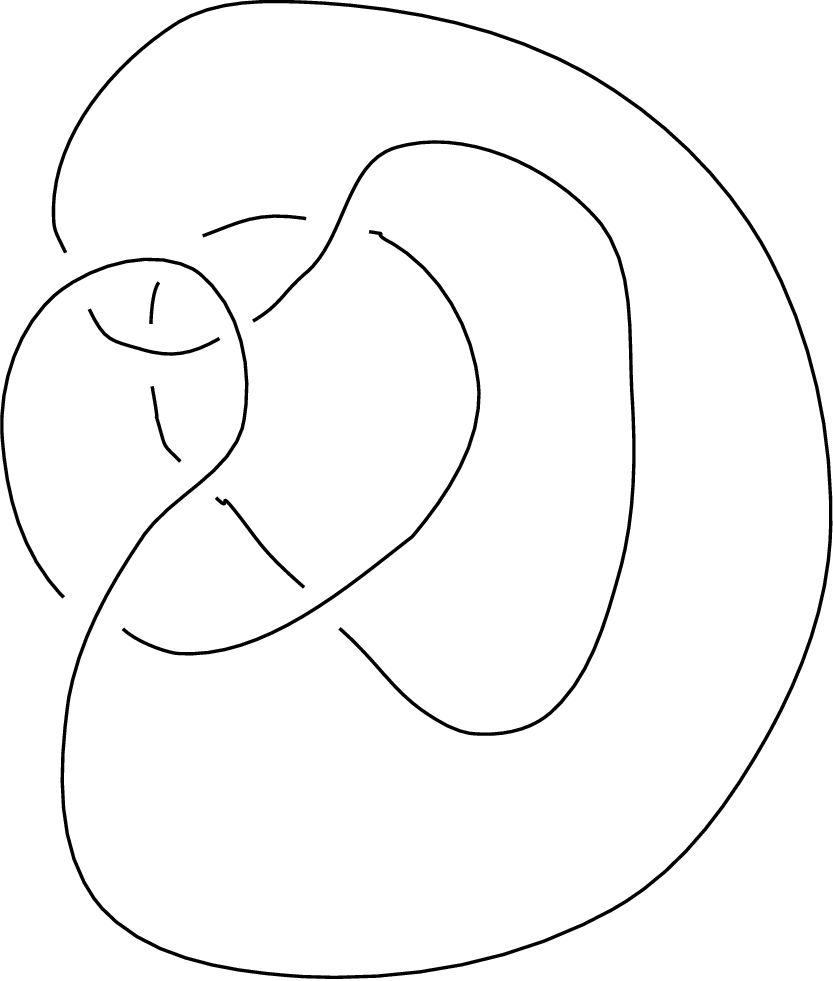}
\caption{A minimal diagram of a trivial knot}
\label{minimal_diagram}
\end{figure}
\end{remark}

\section{Reidemeister move III on the associated minimal diagrams}
\label{sec:relation}

In the above sections, we have considered minimal diagrams by using only Reidemeister moves I and II. Then we may ask what type of Reidemeister move III (and III*) changes the associated minimal diagram.

In this section, by studying Reidemeister moves III and III* from the viewpoint of minimal diagram change, we construct special types of Reidemeister moves III and III* for an arbitrary link without a trivial split component. As a corollary, we may see that, for every RI-II equivalence class (Definition~\ref{def:RI-II equivalence}) of a link without a trivial split component, the set of RI-II equivalence classes $(-)$-adjacent (Definition \ref{def:adjacent}) to the original RI-II equivalence class is the set of RI-II equivalence classes obtained by applying a Reidemeister move III or III* to the minimal diagram in the original RI-II equivalence class.

\begin{definition}
\label{def:RI-II equivalence}
Two link diagrams are said to be \emph{RI-II equivalent} if the diagrams may be transformed to each other by applying finitely many moves RI and RII, without using moves RIII or RIII*.
\end{definition}

\begin{definition}\label{def:minimal_crossing_number}
The number of crossings in the minimal diagrams in a RI-II equivalence class $\mathcal{D}$ is called \emph{crossing number} of $\mathcal{D}$ in this paper and is denoted by $c(\mathcal{D})$ (which is unique because of Corollary \ref{cor:minimal_crossing_number}).
\end{definition}

\begin{definition}\label{def:adjacent}
Two RI-II equivalence classes are defined to be \emph{adjacent}, if a link diagram exists in each of the two RI-II equivalence classes such that they are transformed to each other by applying a single move RIII or RIII*. Moreover, a RI-II equivalence class $\mathcal{D}_{1}$ is defined to be
$(-)$-adjacent (or $(+)$-adjacent) to a RI-II equivalence class $\mathcal{D}_{2}$, if they are adjacent and satisfy $c(\mathcal{D}_{1}) \leq c(\mathcal{D}_{2})$ (resp. $c(\mathcal{D}_{1}) \geq c(\mathcal{D}_{2})$).
\end{definition}

\begin{figure}[H]
\centering
\begin{subfigure}{0.8\linewidth}
\centering
\includegraphics[width=\linewidth]{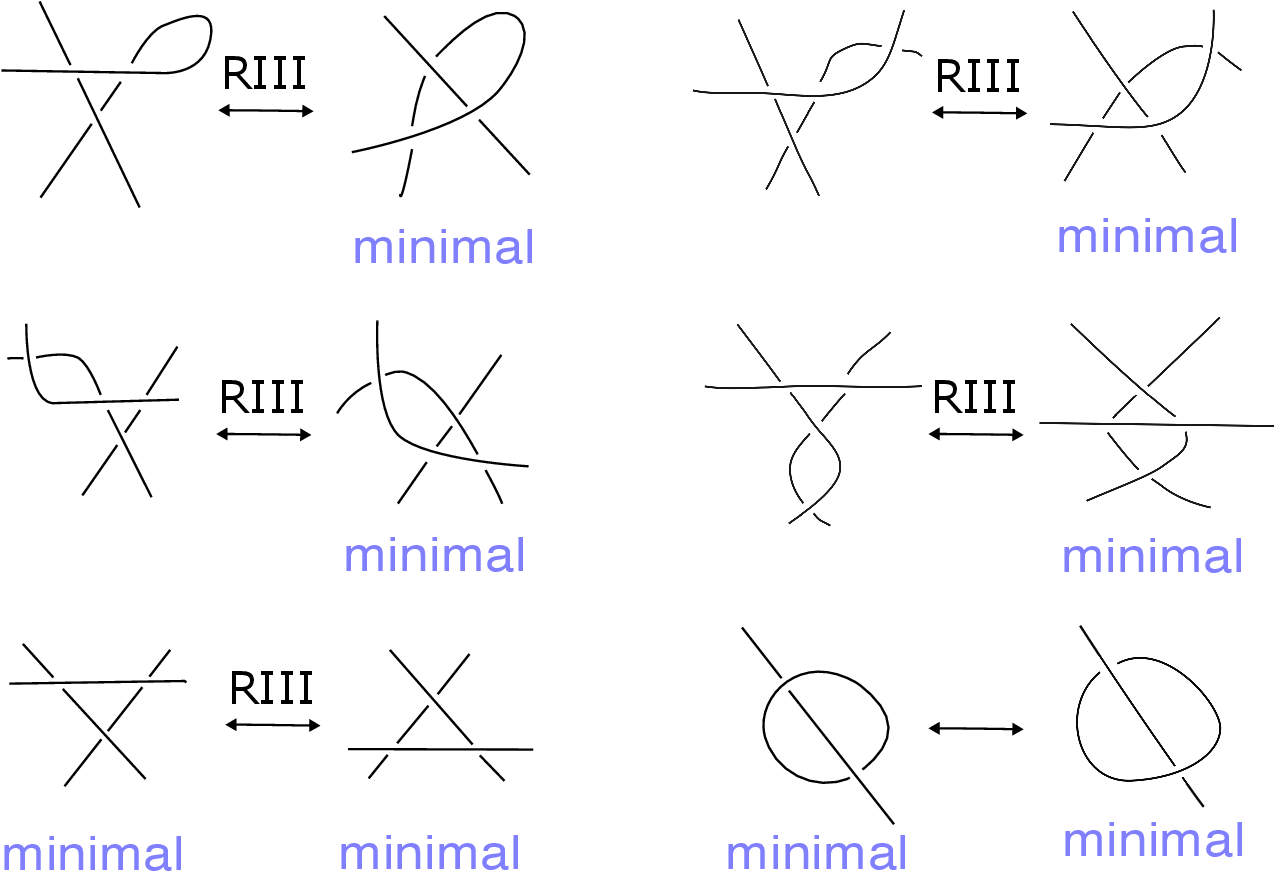}
\caption{RIII's that may change the associated RI-II equivalence classes}
\label{fig:RIII_reduced}
\end{subfigure}
\par\bigskip
\begin{subfigure}{0.8\linewidth}
\centering
\includegraphics[width=\linewidth]{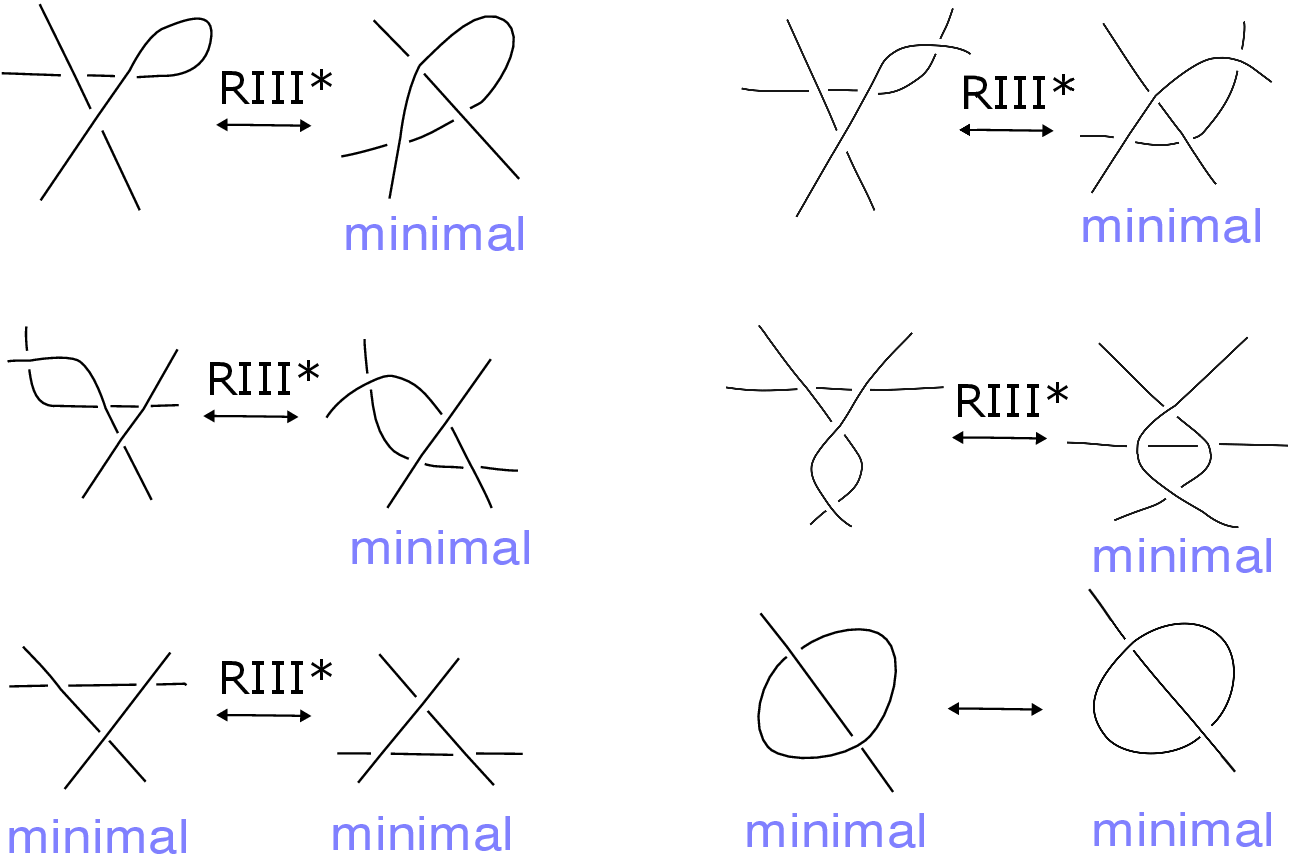}
\caption{RIII*'s that may change the associated RI-II equivalence classes}
\label{fig:RIII*_reduced}
\end{subfigure}
\caption{}
\end{figure}

\begin{theorem}
\label{theo:RIII_reduced}
Assume that two diagrams $D_1$ and $D_2$ of a link without a trivial split component are transformed to each other by applying a single move RIII or RIII*. Then there exist two link diagrams $D_{1}'$ and $D_{2}'$ which are RI-II equivalent to $D_{1}$ and $D_{2}$, respectively, such that they are transformed to each other by applying exactly one of the local moves in Figures~\ref{fig:RIII_reduced}~and~\ref{fig:RIII*_reduced}.
\end{theorem}

\begin{remark}
The diagrams which are indicated as minimal in Figures~\ref{fig:RIII_reduced}~and~\ref{fig:RIII*_reduced} are minimal (see Definition~\ref{def:minimal}). Note that this minimality is not local, but global.
\end{remark}

\begin{remark}\label{rem:RIII_RIII*_reduced}
Every RI-II equivalence class containing the diagram on the left hand side in Figures~\ref{fig:RIII_reduced}~and~\ref{fig:RIII*_reduced} has the same or fewer minimal crossing number (see Definition~\ref{def:minimal_crossing_number}) than the RI-II equivalence class containing the corresponding diagram on the right hand side, since the diagrams on the right hand sides are all minimal and every local move in Figure \ref{fig:RIII_reduced} and \ref{fig:RIII*_reduced} does not change the number of crossings. This means that, if the crossing number of the RI-II equivalence class of $D_{1}$ is equal to or greater than the crossing number of the RI-II equivalence class of $D_{2}$, then $D_{1}$ can be transformed into $D_{2}$ by applying moves $RI_{-}$ and $RII_{-}$ (the diagram becomes minimal), followed by a single move RIII or RIII*, followed by moves $RI_{+}$ and $RII_{+}$.
\end{remark}

\begin{proof}[Proof of Theorem~\ref{theo:RIII_reduced}]
Let us first prove cases of the move RIII. Assume that the two diagrams $D_{1}$ and $D_{2}$ are transformed to each other by applying a single move RIII, not a move RIII*.
When we consider applying finitely many of the same moves $\mathrm{RI}_-$ and $\mathrm{RII}_-$ to the two diagrams $D_{1}$ and $D_{2}$ which change the corresponding positions outside the two local disks of the move RIII until a move $\mathrm{RI}$ or $\mathrm{RII}$ cannot be applied to the outside of the two local disks, we get a special move RIII. Note that only the regions adjacent to the two triangles of the special move RIII can be monogons or the digons for moves RII. Figure~\ref{fig:RIII_special} depicts this special move RIII, whose every region around the two triangles is indicated by one alphabet.

\begin{figure}[H]
\centering
\includegraphics[width=.6\linewidth]{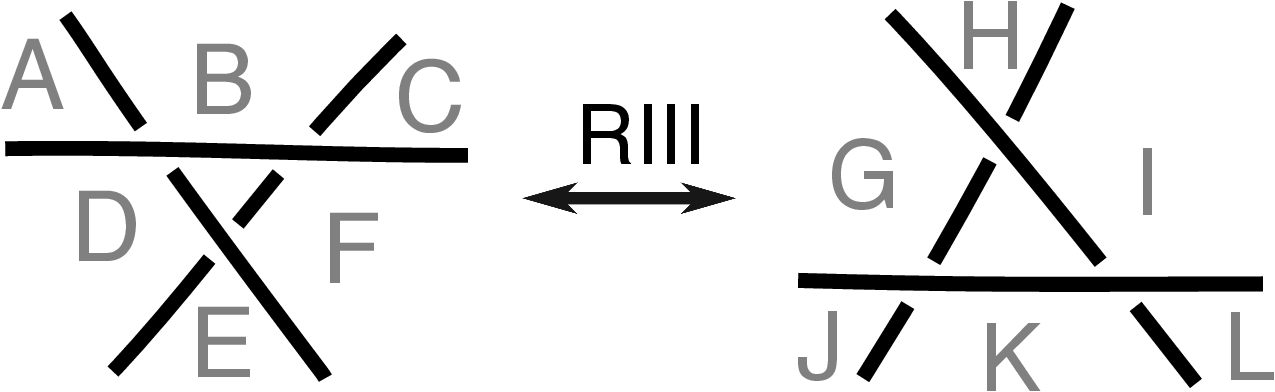}
\caption{The special move RIII where only the regions around the two triangles can be monogons or the digons}
\label{fig:RIII_special}
\end{figure}

Let us classify all the regions adjacent to the two triangles (for the special move RIII depicted in Figure~\ref{fig:RIII_special}) by the numbers of monogons and the digons (for moves RII) in the regions. Assume now that every monogon and special digon of regions $A$--$L$ does cover a single region in the two local disks depicted in Figure~\ref{fig:RIII_special}.

\bigskip

\textbf{(A) The case of the regions adjacent to the two triangles being not monogons or the digons}

In (A), the diagrams on the left and on the right hand sides in Figure~\ref{fig:RIII_special} are both minimal, so the special move RIII is the move of the lower left in Figure~\ref{fig:RIII_reduced}.

\medskip

\textbf{(B) The cases of the regions adjacent to the two triangles being not the digons}

In (B), the region which can be a monogon is only region $C$ or $J$. 
This reason is explained right below.

First, every region in regions $B$, $D$, $F$, $G$, $I$ and $K$ contains two or more crossings, so these regions cannot be monogons.

Second, when one of regions $A$, $E$, $H$, $L$ is a monogon, automatically a region becomes the digon. For instance, when region $A$ is a monogon, automatically region $G$ becomes the digon. So, these cases are not included in (B).

Let us classify (B) by the number of regions $A$--$L$ being monogons (The number is one or two, because only region C or J can be a monogon in (B)).

\medskip

\uline{(B-1) When the number of monogons in the regions adjacent to the two triangles is only one.}

When only region $C$ is a monogon, the diagram on the right hand side in Figure \ref{fig:RIII_special} is minimal. Hence, in this case, the special move RIII is the move of the upper left in Figure~\ref{fig:RIII_reduced}.

When only region $J$ is a monogon, the diagram on the left hand side in Figure \ref{fig:RIII_special} is minimal. Hence, in this case, the special move RIII is the move of the upper left in Figure~\ref{fig:RIII_reduced}.

\medskip

\uline{(B-2) When the number of monogons in the regions adjacent to the two triangles is two.}

When only regions $C$ and $J$ are both monogons, as in Figure~\ref{fig:RIII_example1}, the special move RIII can be expressed by applying finitely many moves RI and the move of the lower right in Figure \ref{fig:RIII_reduced}. Note that the diagrams at the lower left and the lower right in Figure~\ref{fig:RIII_example1} are both minimal, since the regions which can be seen in the two local disks cannot be monogons or the digons.

\begin{figure}[H]
\centering 
\includegraphics[width=.7\linewidth]{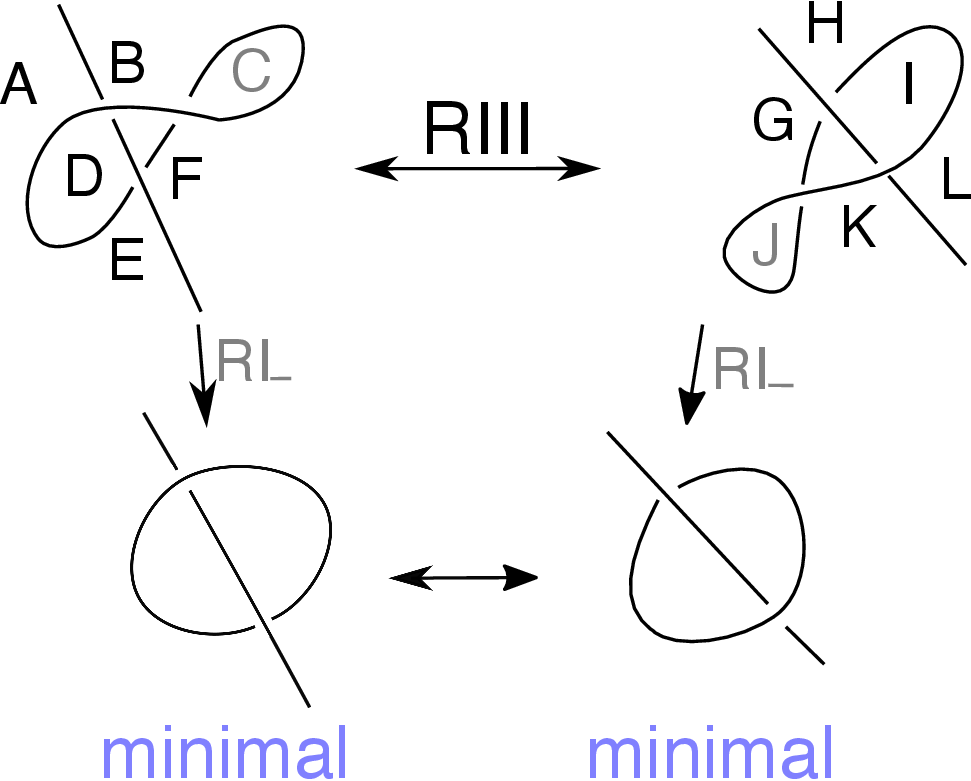}
\caption{Regions $C$ and $J$ are monogons}
\label{fig:RIII_example1}
\end{figure}

\textbf{(C) The cases of the regions adjacent to the two triangles being not monogons} 

In (C), the region which can be the digon is region $A$, $C$, $E$, $H$, $J$ or $L$ in Figure \ref{fig:RIII_special}. 
This reason is explained right below.

First, region $D$ or $I$ cannot be the digon, because of the crossings appeared in the local disks. 

Second, when region $B$, $F$, $G$ or $K$ is the digon, a region automatically becomes a monogon. For instance, when region $B$ is the digon, region $H$ automatically becomes a monogon. So, these cases are not included in (C).  

Let us classify (C) by the number of regions $A$--$J$ being the digons.

\medskip

\uline{(C-1) When the number of the digons in the regions adjacent to the two triangles is only one.}

When only region $A$ is the digon, the diagram on the right hand side in Figure \ref{fig:RIII_special} is minimal. Hence, the special move RIII is the move of the middle left in Figure~\ref{fig:RIII_reduced}.

The proofs of when only regions $C$, $E$, $H$, $J$ and $L$ are the digons are clearly the same as the right above proof. 

\medskip

\uline{(C-2) When the number of the digons in the regions adjacent to the two triangles is two.}

The first cases that we consider in (C-2) are the pairs $A$$C$, $A$$E$, $C$$E$, $H$$J$, $J$$L$ and $H$$L$ being the digons ($6$ cases). Note that, only in these cases, the diagrams on the other hand side in Figure \ref{fig:RIII_special} are all minimal.

The second cases that we consider in (C-2) are the pairs $A$$H$, $A$$J$, $A$$L$, $C$$H$, $C$$J$, $C$$L$, $E$$H$, $E$$J$ and $E$$L$ being the digons ($9$ cases). Note that, only in these cases, the diagrams on the both hand sides in Figure \ref{fig:RIII_special} are not minimal. 

Let us start the first cases. When only regions $A$$C$ are the digons, the diagram on the right hand side in Figure \ref{fig:RIII_special} is minimal, so the special move RIII is the move of the upper right or the middle left in Figure~\ref{fig:RIII_reduced}.
The proofs of when pairs $A$$E$, $C$$E$, $H$$J$, $J$$L$ and $H$$L$ are the digons are clearly the same as the right above proof and the special moves RIII are the moves of the upper right and the middle left and the middle right in Figure~\ref{fig:RIII_reduced}. We have checked all the $6$ cases in the first cases of (C-2).

Let us now check the second cases of (C-2). When regions $A$$H$ or $E$$L$ are the digons, the two cases cannot occur. See Figure~\ref{fig:RIII_example2}~and~Figure~\ref{fig:RIII_proof_example5}, respectively. In Figure~\ref{fig:RIII_example2}, region $H$ cannot be the digon when region $A$ is the digon. In Figure \ref{fig:RIII_proof_example5}, region L cannot be the digon when region E is the digon.

\begin{figure}[H]
\centering
\includegraphics[width=.6\linewidth]{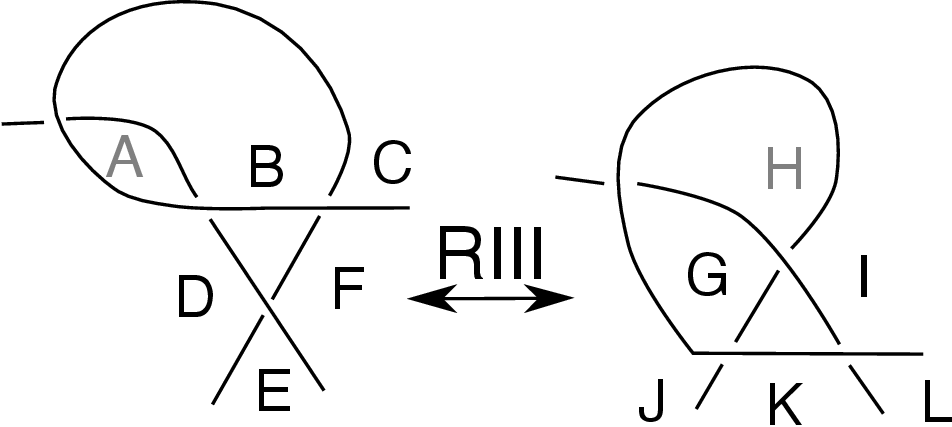}
\caption{When region $A$ is the digon and region $H$ has two crossings}
\label{fig:RIII_example2}
\end{figure}

\begin{figure}[H]
\centering
\includegraphics[width=.6\linewidth]{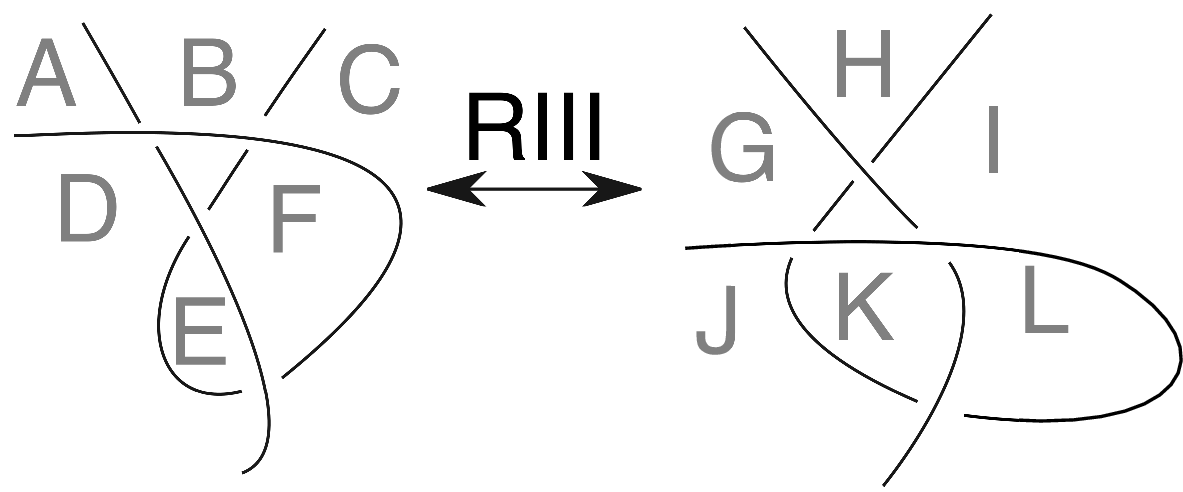}
\caption{When region $E$ is the digon and region $L$ has two crossings}
\label{fig:RIII_proof_example5}
\end{figure}

When pairs $A$$J$, $C$$H$, $C$$L$ and $E$$J$ are the digons, the special moves RIII can be expressed by applying finitely many moves RI and RII. See Figures~\ref{fig:RIII_proof_example8}, \ref{fig:RIII_proof_example7}, \ref{fig:RIII_proof_example3}~and~\ref{fig:RIII_proof_example6}, respectively. So, these cases are not included in Theorem~\ref{theo:RIII_reduced}.

\begin{figure}[H]
\centering
\begin{minipage}[t]{0.49\textwidth}
\centering
\includegraphics[width=0.8\linewidth, height=6cm]{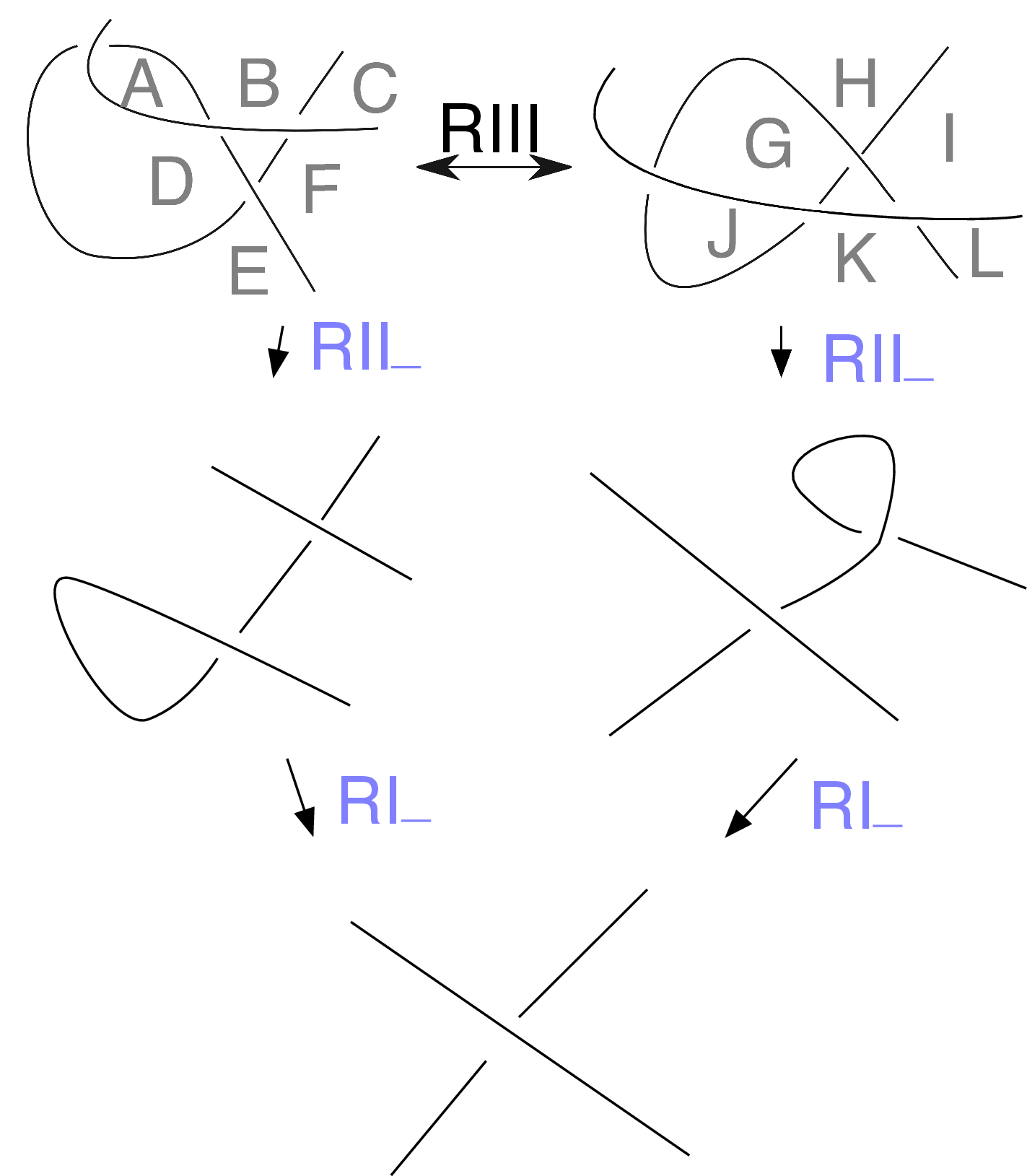}
\caption{Regions $A$$J$ are digons}
\label{fig:RIII_proof_example8}
\end{minipage}
\begin{minipage}[t]{0.49\textwidth}
\centering
\includegraphics[width=0.8\linewidth, height=6cm]{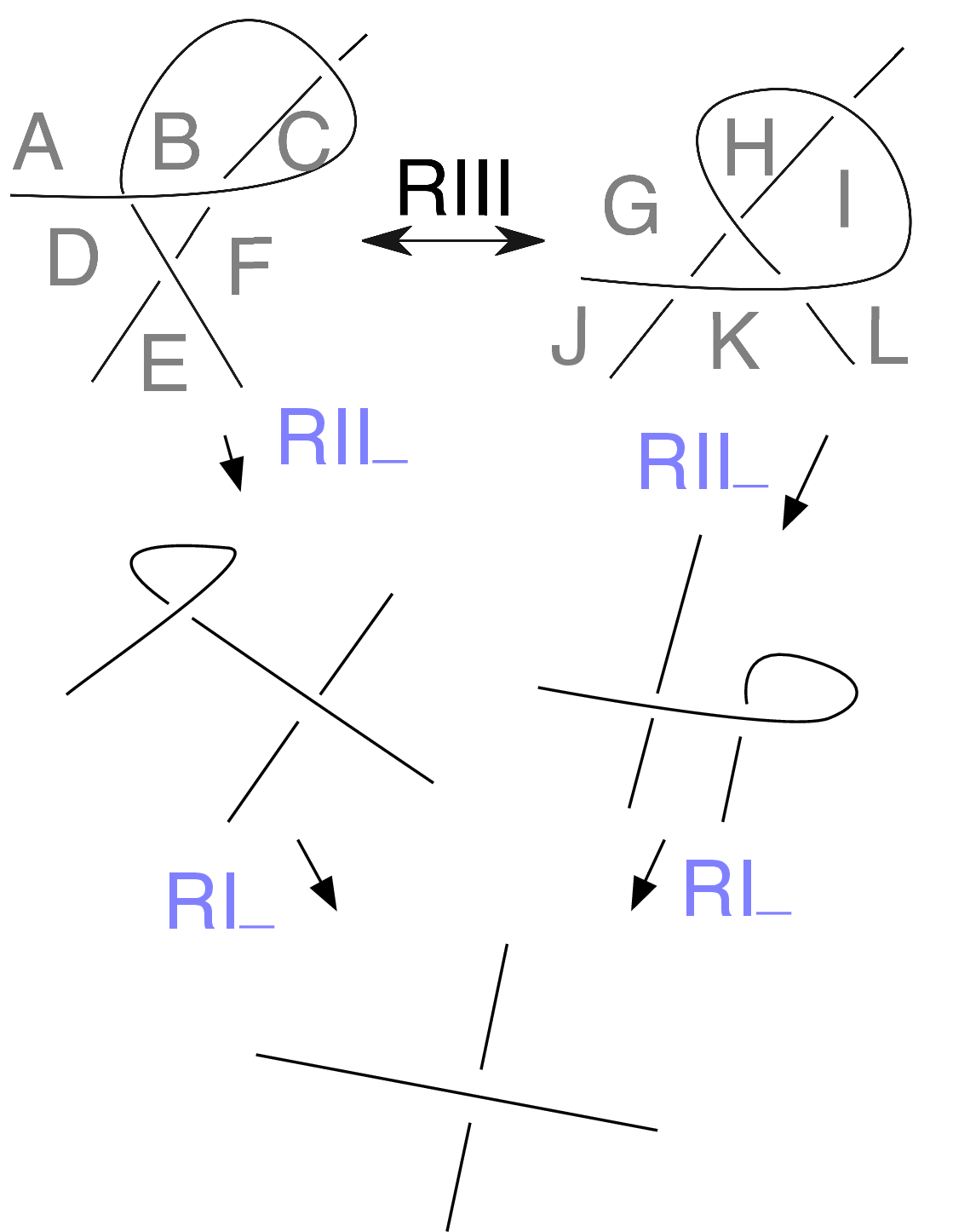}
\caption{Regions $C$$H$ are the digons}
\label{fig:RIII_proof_example7}
\end{minipage}
\end{figure}

\begin{figure}[H]
\centering
\begin{minipage}[t]{0.49\textwidth}
\centering
\includegraphics[width=.8\linewidth, height=6cm]{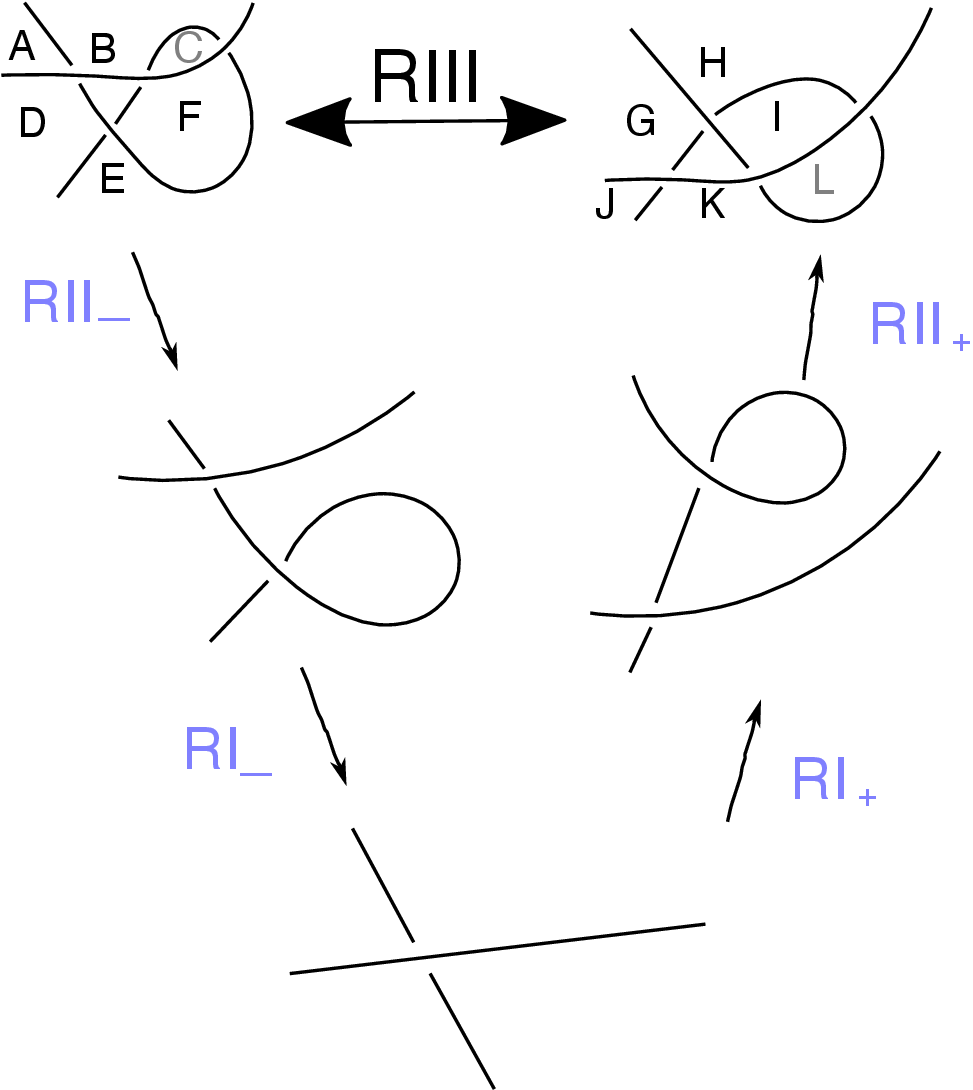}
\caption{Regions $C$$L$ are the digons}
\label{fig:RIII_proof_example3}
\end{minipage}
\begin{minipage}[t]{0.49\textwidth}
\centering
\includegraphics[width=.8\linewidth, height=6cm]{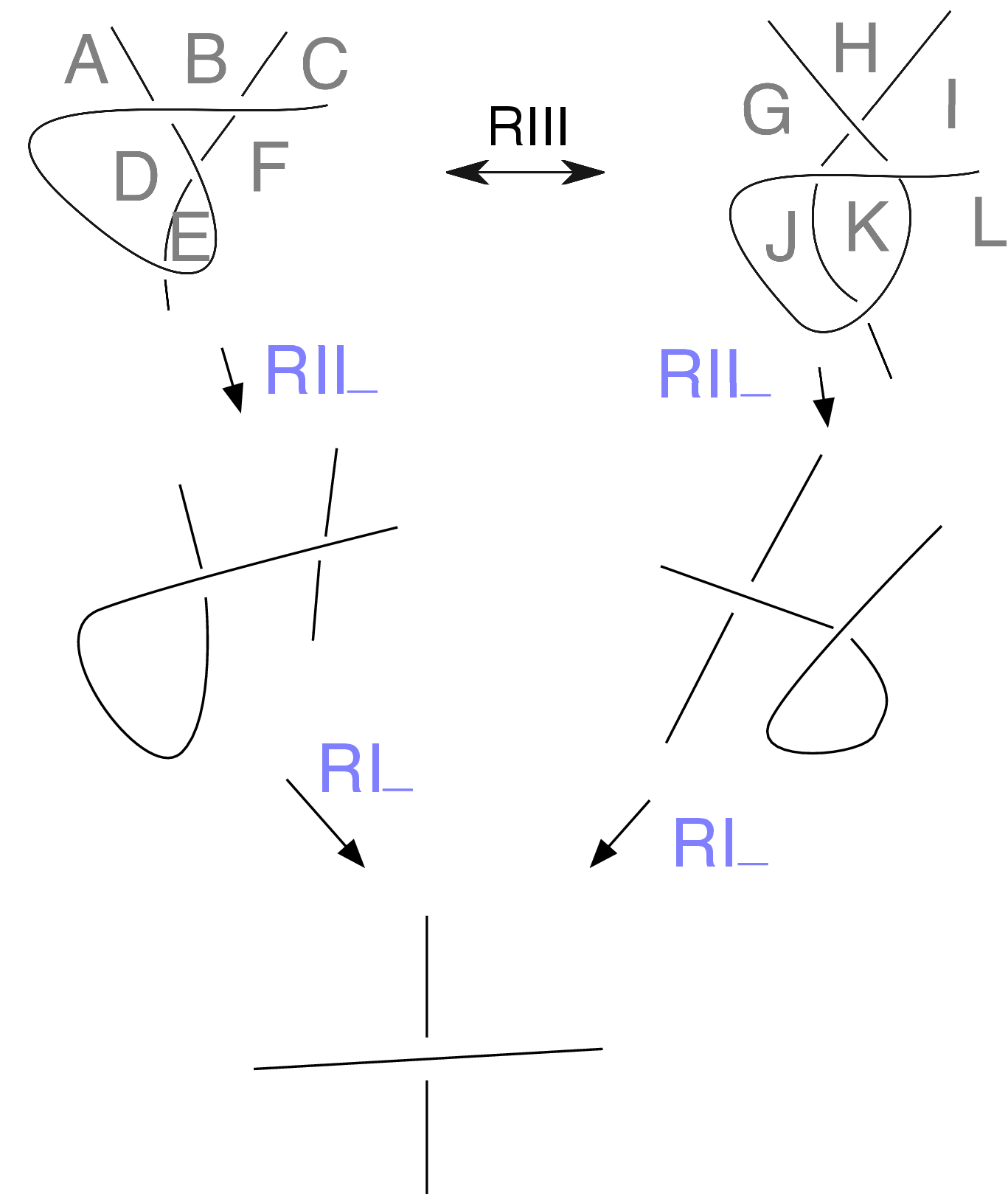}
\caption{Regions $E$$J$ are the digons}
\label{fig:RIII_proof_example6}
\end{minipage}
\end{figure}

When pairs $A$$L$, $C$$J$ and $E$$H$ are the digons, the special moves RIII can be expressed by applying finitely many moves RII and a single move RIII*. See Figure~\ref{fig:RIII_proof_example9}, \ref{fig:RIII_proof_example10},~and~\ref{fig:RIII_proof_example11}, respectively. So, these cases result in cases of the move RIII*. Note that, these processes which change moves RIII into moves RIII* always decrease the number of crossings. We deal with this kind of processes after we discuss the other cases of the move RIII.

\begin{figure}[H]
\centering
\begin{minipage}[t]{0.49\textwidth}
\centering
\includegraphics[width=.8\linewidth, height=5cm]{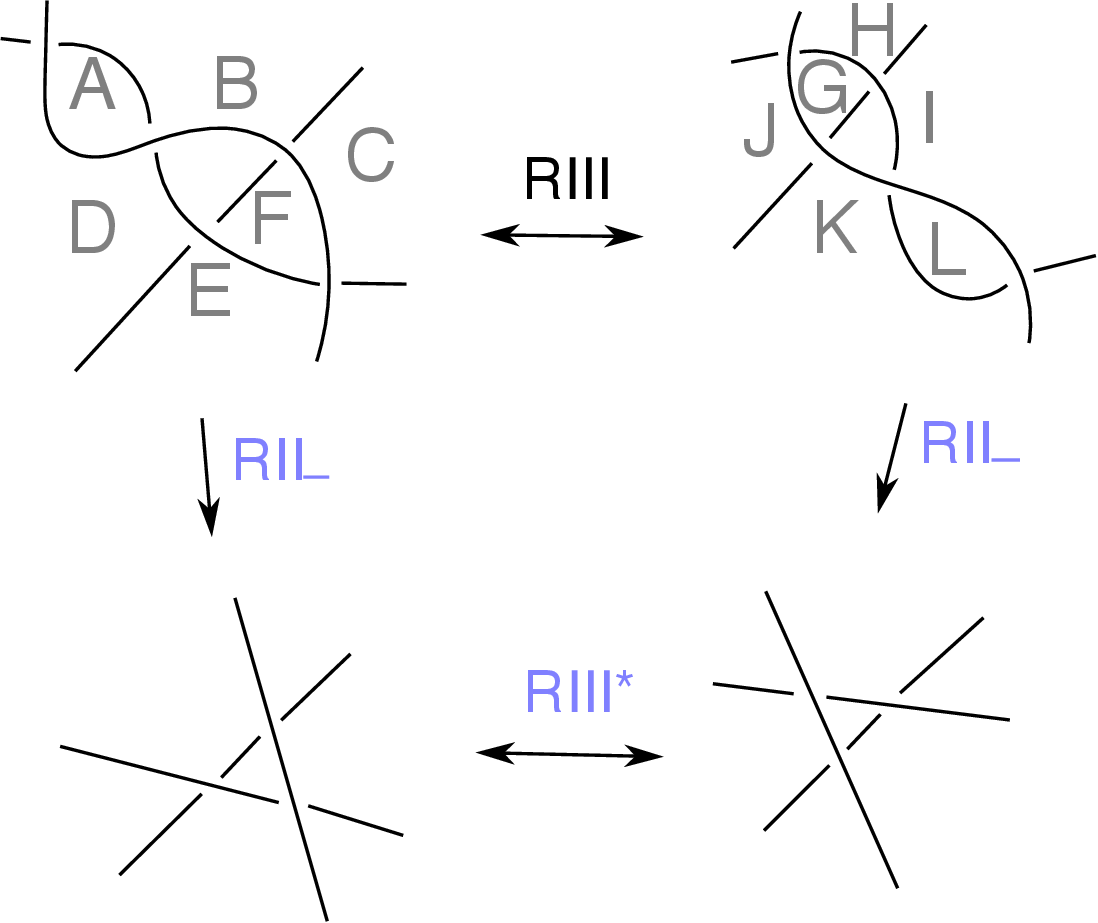}
\caption{Regions $A$$L$ are the digons}
\label{fig:RIII_proof_example9}
\end{minipage}
\begin{minipage}[t]{0.49\textwidth}
\centering
\includegraphics[width=.8\linewidth, height=5cm]{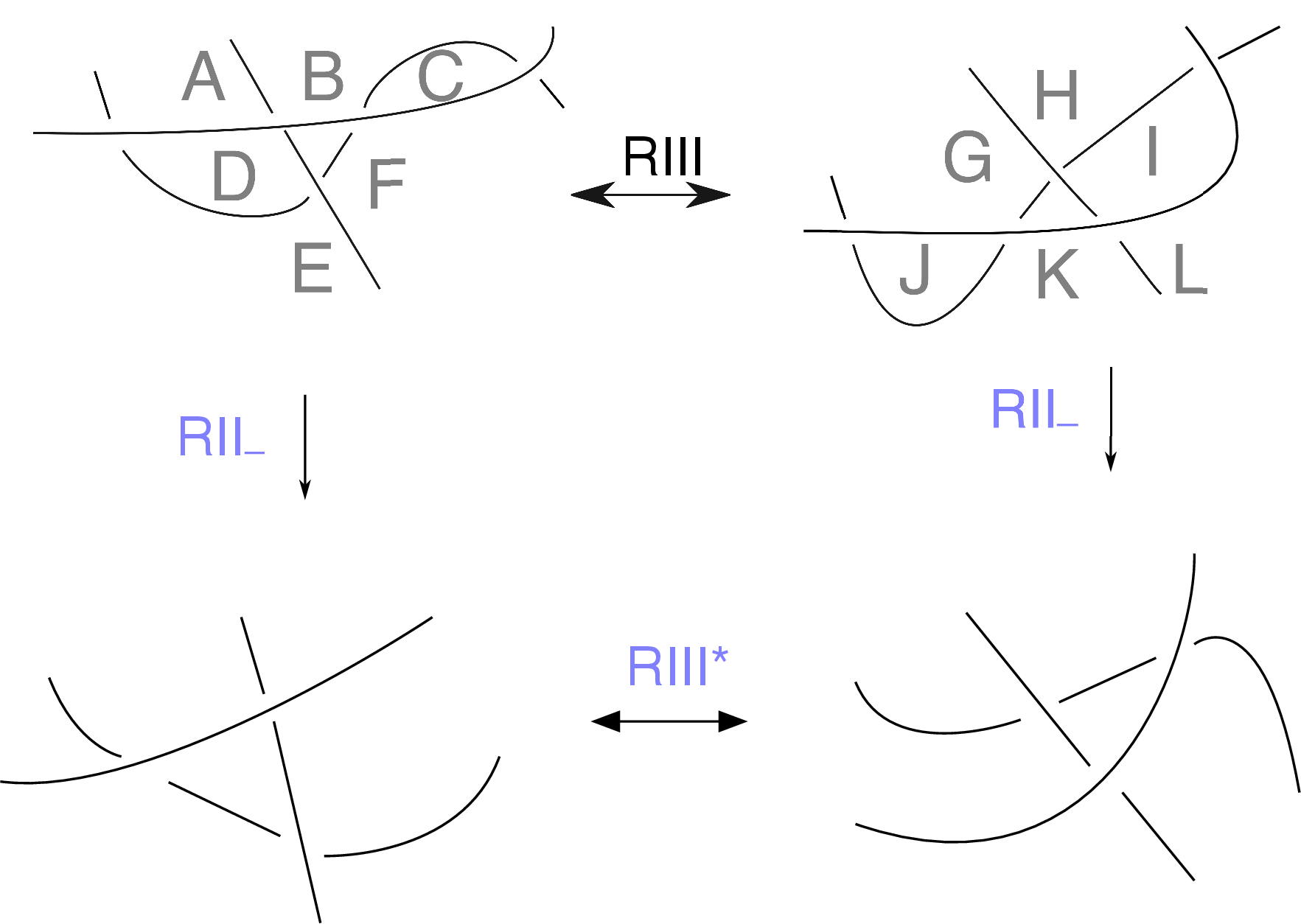}
\caption{Regions $C$$J$ are the digons}
\label{fig:RIII_proof_example10}
\end{minipage}
\end{figure}

\begin{figure}[H]
\centering
\includegraphics[width=.5\linewidth,keepaspectratio]{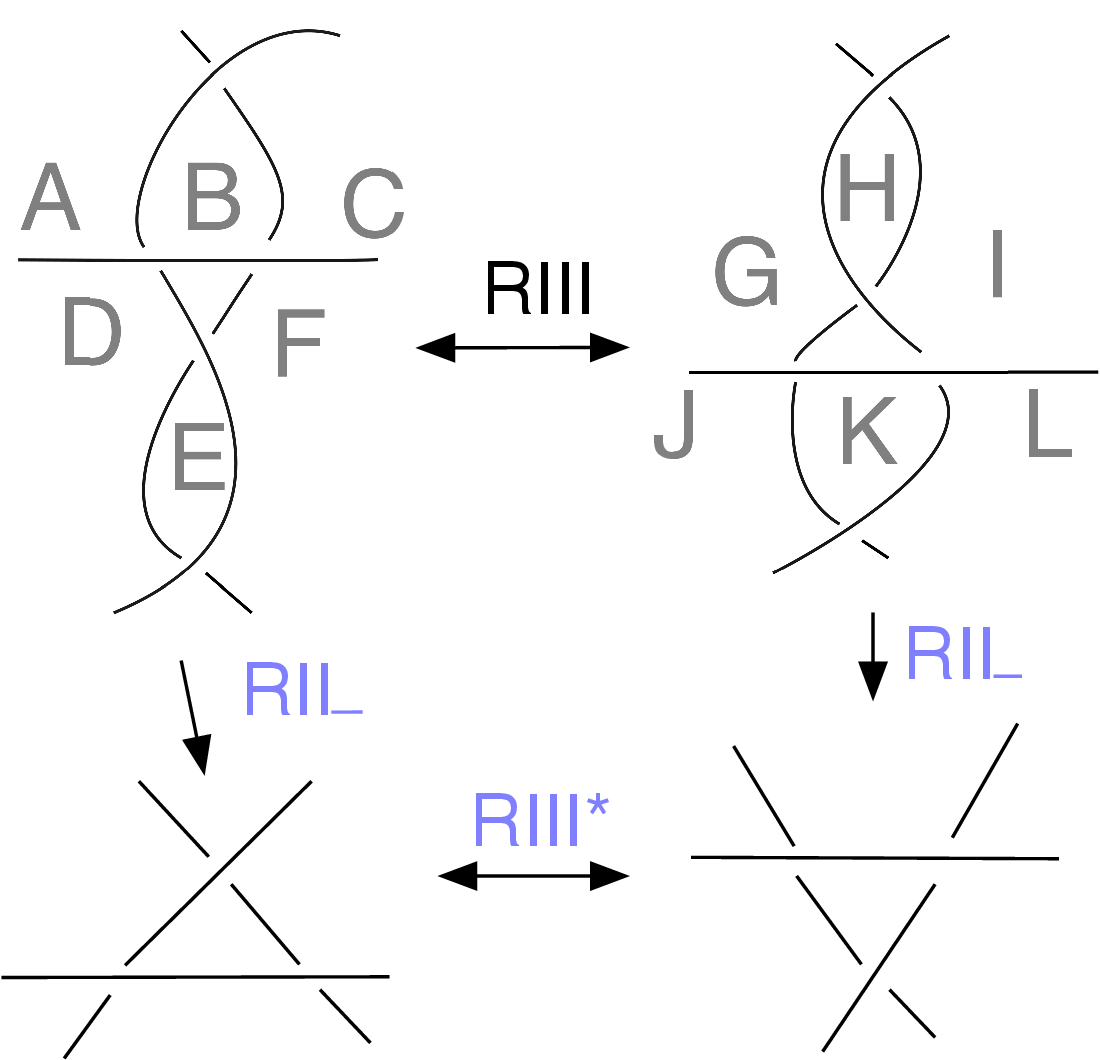}
\caption{Regions $E$$H$ are the digons}
\label{fig:RIII_proof_example11}
\end{figure}

We have checked all the $9$ second cases of (C-2) right above. 

All the cases in (C-2) have been checked. 

\medskip

\uline{(C-3) When the number of the regions adjacent to the two triangles being the digons is three or more.}

In (C-3), we consider combinations of regions A-L being the digons.

The first cases which we consider in (C-3) are the combinations existing only on the diagram of the left or the right hand side in Figure \ref{fig:RIII_special} (in other words, the two triples $A$$C$$E$ and $H$$J$$L$ being the digons). Note that the diagram on the other hand side is minimal.

The second cases which we consider in (C-3) are the combinations existing on the diagrams of both the left and the right hand sides. We can recall here that in the second cases of (C-2), we used only the information of the positions of the digons, which exist on both the left and the right hand sides in Figure \ref{fig:RIII_special}. Hence, the second cases here result in (C-2).

When the triple $A$$C$$E$ are the digons, the diagram on the right hand side is minimal. So the special move RIII is the move of the top right or the middle left or the middle right in Figure~\ref{fig:RIII_reduced}.

The case of the triple $H$$J$$L$ being the digons is clearly the same as the right above case. The special move RIII is the move of the top right or the middle left or the middle right in Figure~\ref{fig:RIII_reduced}.

All the cases in (C-3) have been checked.

\bigskip

\textbf{(D) The cases of there existing both monogons and the digons in the regions adjacent to the two triangles.}

The region which may be a monogon in Figure \ref{fig:RIII_special} is region $A$, $C$, $E$, $H$, $J$ or $L$. We classify (D) by which regions among regions $A$, $C$, $E$, $H$, $J$ and $L$ being monogons.

When region $A$ is a monogon (automatically region $G$ becomes the digon), as in Figure~\ref{fig:RIII_proof_example4}, the special move RIII may be expressed by applying finitely many moves RI and RII. So, this case is not included in Theorem~\ref{theo:RIII_reduced}. 

\begin{figure}[H]
\centering
\includegraphics[width=.7\linewidth]{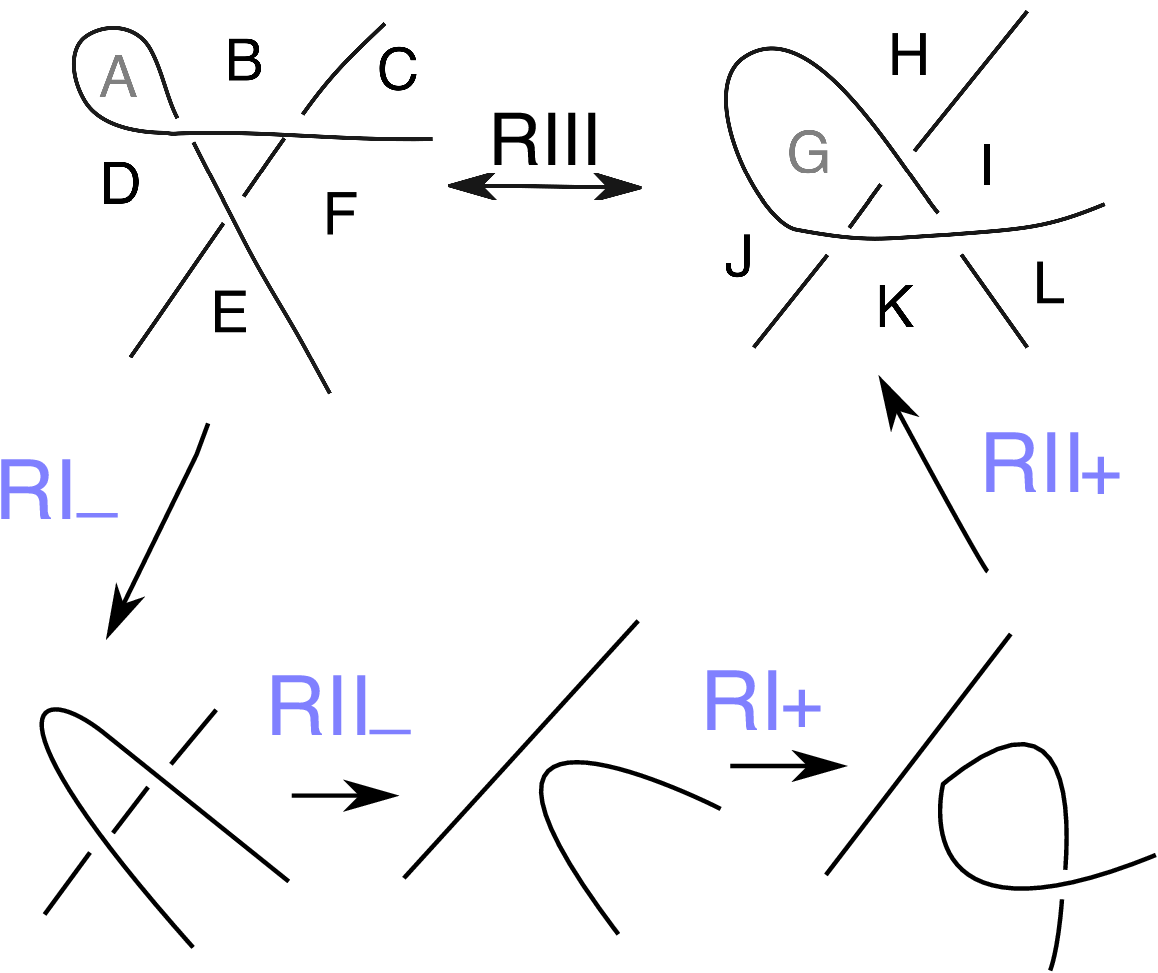}
\caption{Region $A$ is a monogon (automatically, region $G$ becomes the digon)}
\label{fig:RIII_proof_example4}
\end{figure}

The proofs of regions $E$, $H$ and $L$ being monogons are the same as the right above proof and the special moves RIII can be expressed by applying finitely many moves RI and RII. Hence, these cases are not included in Theorem~\ref{theo:RIII_reduced}.

By the right above proofs, we can disregard when regions $A$, $E$, $H$ and $L$ are monogons. Hence, we regard only when region $C$ or $J$ is a monogon from here.

\medskip

\uline{(D-1) The cases that region C is a monogon and another region in Figure \ref{fig:RIII_special} is not a monogon.}

Assume now that region $C$ is a monogon and another region in Figure~\ref{fig:RIII_special} is not a monogon. Regions $H$, $L$, $B$ and $F$ in Figure~\ref{fig:RIII_special} cannot be the digons, since the assumption of region $C$ being a monogon affects the regions around region $C$, and region $D$ cannot be the digon because of the crossings appeared in the local disk. When region $G$ or $K$ is the digon, the special move RIII can be expressed by applying finitely many moves RI and RII. See Figure~\ref{fig:RIII_proof_example4}. So, the region which can be the digon in (D-1) is region $A$, $E$ or $J$. Let us consider the combinations of the regions $A$, $E$ and $J$ being the digons (in other words, $A$, $E$, $J$, $A$$E$, $A$$J$, $E$$J$ and $A$$E$$J$ being the digons ($7$ cases)). We can recall that, when combinations $A$$J$, $E$$J$ and $A$$E$$J$ are the digons, the special moves RIII can be expressed by applying finitely many moves RI and RII. See Figures~\ref{fig:RIII_proof_example8}~and~\ref{fig:RIII_proof_example6}. So it is enough to check the cases that combinations $A$, $E$, $J$ and $A$$E$ are the digons. When combinations $A$, $E$ and $A$$E$ are the digon(s), the diagrams on the right hand sides in Figure \ref{fig:RIII_special} are all minimal, so the special moves RIII are the moves of the top left and the middle left and the middle right in Figure~\ref{fig:RIII_reduced}. When region $J$ is the digon, the special move RIII can be expressed by applying a single move RIII* and finitely many moves RI and RII. See Figure~\ref{fig:RIII_proof_example12}. This case result in the cases of the move RIII*. Note that, this process which change the move RIII into the move RIII* decreases the number of crossings.

\begin{figure}[H]
\centering
\includegraphics[width=.6\linewidth,keepaspectratio]{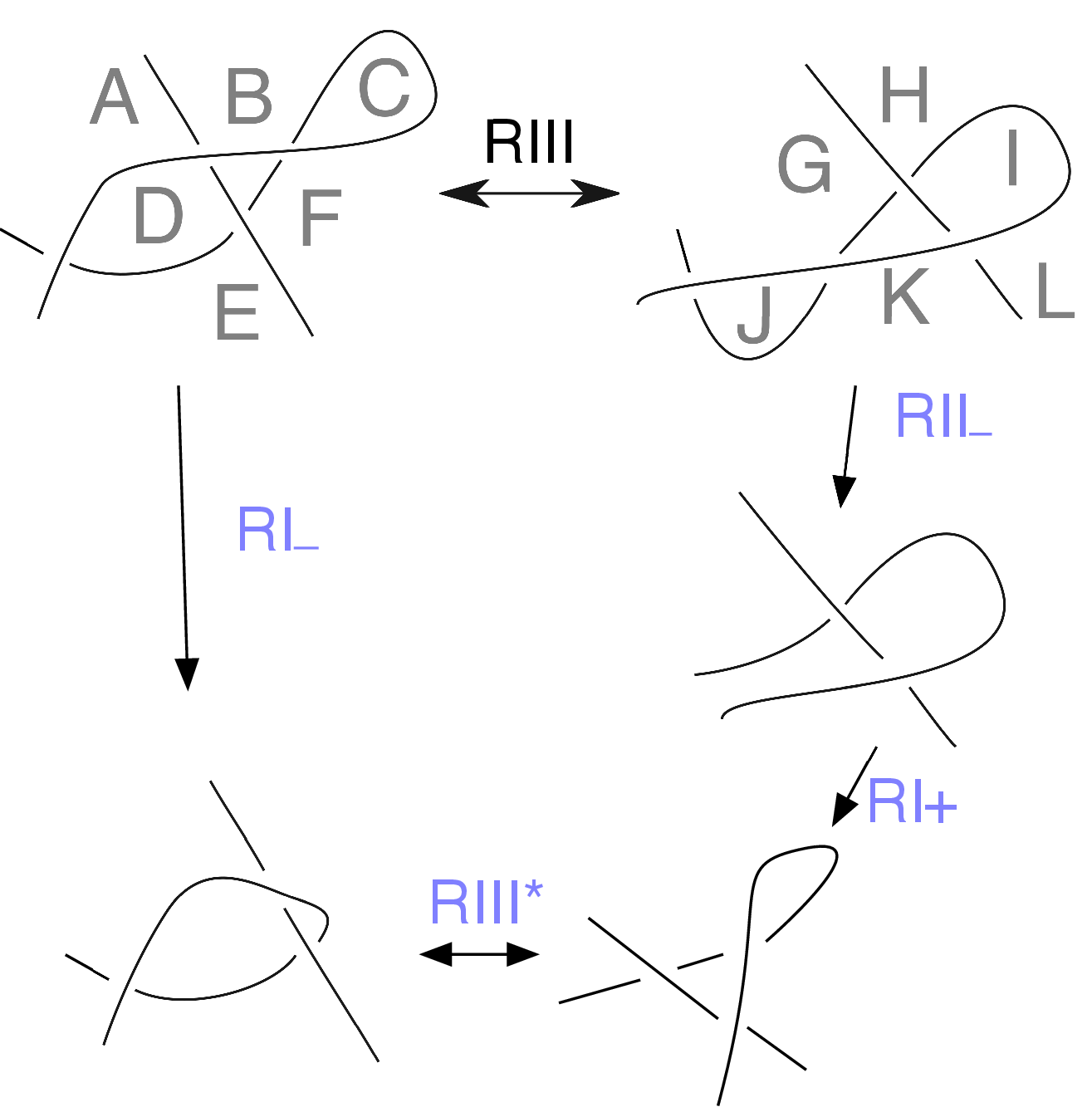}
\caption{Region $C$ is a monogon, region $J$ is the bigon}
\label{fig:RIII_proof_example12}
\end{figure}

\medskip

\uline{(D-2) The cases that region $J$ is a monogon and no another region in Figure \ref{fig:RIII_special} is a monogon}

The proof of region $J$ being a monogon and no another region being a monogon in Figure \ref{fig:RIII_special} is the same as the right above proof, since the position of region $J$ is the same as the one of region $C$ for the special move RIII.
 
 \medskip
 
\uline{(D-3) The cases that both regions C and J are monogons and no another region is a monogon.}
 
The last cases that we should check in (D) is (D-3). Assume now that regions $C$ and $J$ are monogons and no another region is a monogon in Figure \ref{fig:RIII_special}. In this case, every other region (in other words, region $A$, $B$, $E$, $F$, $G$, $H$, $K$ or $L$) cannot be the digon, since the assumption of regions $C$ and $J$ being monogons affects the regions around regions $C$ and $J$. So this case is not included in (D).

\medskip

(A), (B), (C) and (D) above are all the cases that every monogon and special digon in regions $A$--$L$ does cover a single region in the two local disks in Figure~\ref{fig:RIII_special}. What we should check after that is the cases that monogons and the digons in regions $A$--$J$ that do not cover a single region in the two local disks in Figure~\ref{fig:RIII_special}. A monogon contains one crossing, so a monogon cannot cover two regions in the two local disks. By considering that the digon contains two crossings, what we should check here are the cases that pairs $A$$C$, $C$$E$, $H$$J$ and $J$$L$ are the digons which cover the pairs. (Note here that pairs $A$$E$ and $H$$L$ cannot be the digons which cover the pairs, since the information of the crossings that the pairs $A$$E$ and $H$$L$ contain are different from the ones of the digon.) In fact, the cases cannot occur, since in the cases, trivial split components appear near the local disks. For instance, see Figure~\ref{fig:RIII_proof_example13}, where there might be a link diagram inside the dotted circle. 

\begin{figure}[H]
\centering
\includegraphics[width=.6\linewidth]{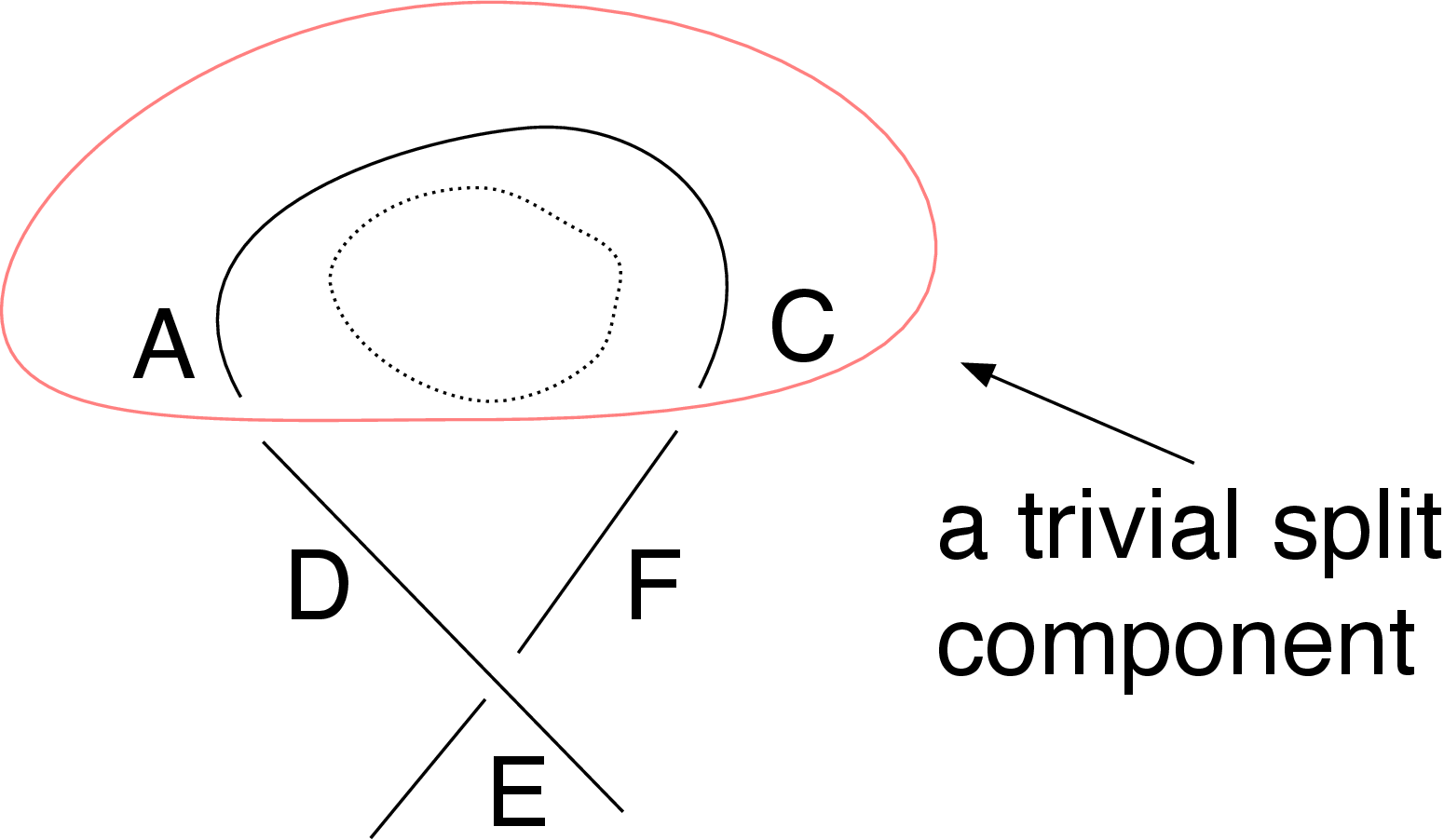}
\caption{Regions $A$ and $C$ are the single digon}
\label{fig:RIII_proof_example13}
\end{figure}

Above are all the cases for the move RIII. Note that the processes appeared above which change moves RIII into moves RIII* always decrease the numbers of crossings. 

\bigskip

Let us now prove the case of the move RIII* as the proof of the move RIII. The move RIII* depicted in Figure~\ref{fig:RIII*_special} is obtained from the move RIII* by applying finitely many of the same moves $\mathrm{RI}_-$, $\mathrm{RII}_-$ to the outside of the two local disks until a move $\mathrm{RI}$ or $\mathrm{RII}$ cannot be applied to the outside of the two local disks, where every region adjacent to the two triangles is indicated by one alphabet. We say that the move RIII* depicted in Figure~\ref{fig:RIII*_special} is "the special move RIII*". Note that Figure~\ref{fig:RIII*_special} is the mirror image of Figure~\ref{fig:RIII_special} including alphabets.

\begin{figure}[H]
\centering
\includegraphics[width=.6\linewidth]{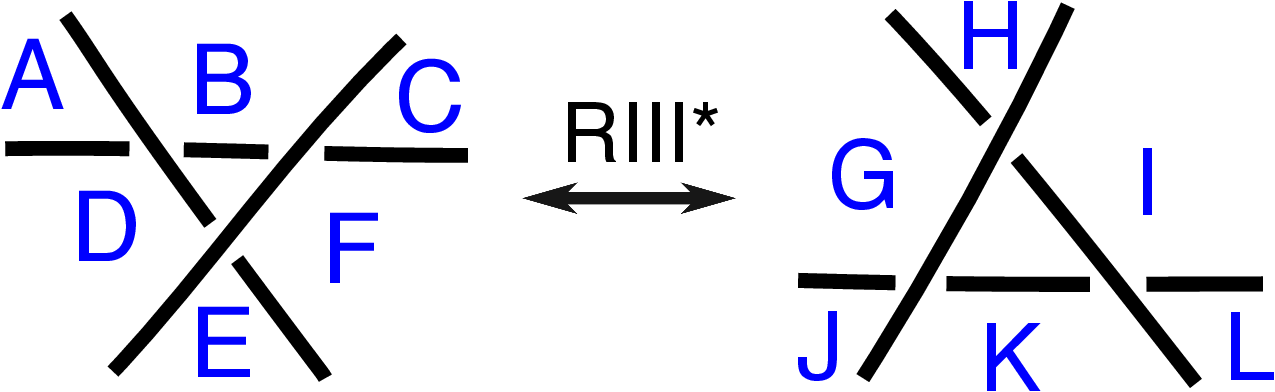}
\caption{The special move RIII* - mirror image of Figure \ref{fig:RIII_special}}
\label{fig:RIII*_special}
\end{figure}  

In fact, every case of the special move RIII* corresponds to a case of the special move RIII. For instance, Figure~\ref{fig:RIII*_proof_example1} is the mirror image of Figure~\ref{fig:RIII_proof_example3}. All the cases of the special move RIII* can be proved as this, since the moves RI, RII, RIII, RIII* also can be applied to the mirror images (The moves RIII change into the moves RIII*, and vice versa). Hence, the result of all the cases of the move RIII* is the mirror image of Figure~\ref{fig:RIII_reduced}, in other words, Figure~\ref{fig:RIII*_reduced}.

The processes which change moves RIII* into moves RIII also appear in the cases of the special move RIII*, and always decrease the number of crossings. The processes which changes moves RIII into moves RIII* and moves RIII* into moves RIII, always decrease the numbers of crossings, so we can see that the cases result in the other cases, since the number of crossings of every link diagram is finite. This means that we can disregard these cases, which completes the proof.

\begin{figure}[H]
\centering
\includegraphics[width=.7\linewidth]{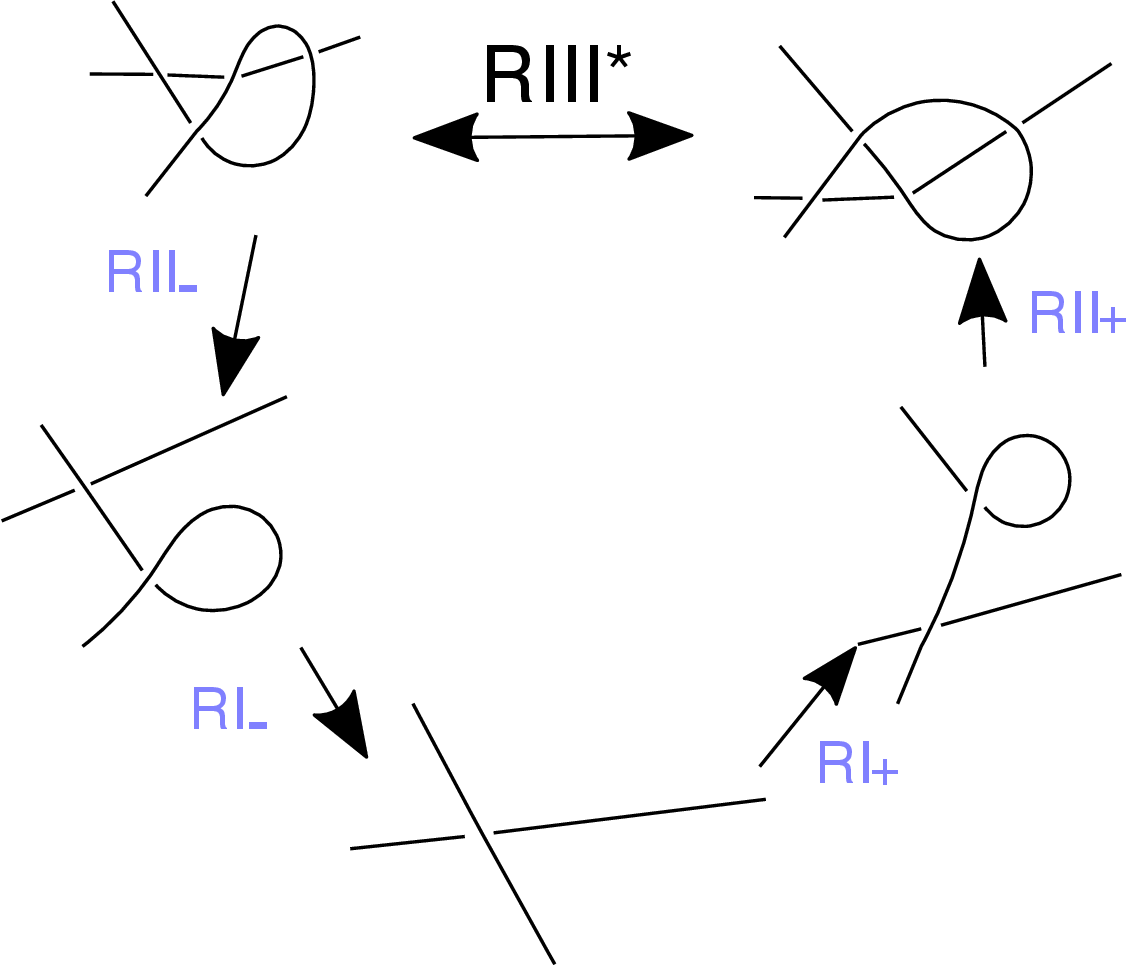}
\caption{Mirror image of Figure \ref{fig:RIII_proof_example3} (Regions CL are the digons)}
\label{fig:RIII*_proof_example1}
\end{figure}
\end{proof}

\begin{corollary}
\label{cor:equivalence_relation}
Let $\mathcal{R}$ be an arbitrary RI-II equivalence class of a link without a trivial split component. We denote the set of all the RI-II equivalence classes $(-)$-adjacent to the RI-II equivalence class $\mathcal{R}$ by $\textbf{A}$($\mathcal{R}$;$-$), and the set of the RI-II equivalence classes containing all the diagrams obtained by applying a Reidemeister move III or III* to the minimal diagram in the RI-II equivalence class $\mathcal{R}$ by $\textbf{M$(\mathcal{R})$}$, then
$\textbf{A}$($\mathcal{R}$;$-$) $=$ $\textbf{M$(\mathcal{R})$}$.
\end{corollary}

\begin{remark}
Due to Theorem \ref{theo:RIII_reduced}, the RI-II equivalence class $\mathcal{R}$ in Corollary \ref{cor:equivalence_relation} has an unique minimal diagram. 
\end{remark}

\begin{proof}
To prove $\textbf{M$(\mathcal{R})$}$ $\subset$ $\textbf{A}$($\mathcal{R}$;$-$) is trivial by the definitions. By Remark \ref{rem:RIII_RIII*_reduced} of Theorem \ref{theo:RIII_reduced}, every RI-II equivalence class in $\textbf{A}$($\mathcal{R}$;$-$) is contained in $\textbf{M$(\mathcal{R})$}$, which proves $\textbf{A}$($\mathcal{R}$;$-$) $\subset$ $\textbf{M$(\mathcal{R})$}$.   
\end{proof}


\begin{thebibliography}{99}

\bibitem{1}A.~Coward, M.~Lackenby: An upper bound on Reidemeister moves, Amer. J. Math. 136 (2014), no.~4, 1023--1066.

\bibitem{2}C.~Hayashi: The number of Reidemeister moves for splitting a link, Math. Ann. 332 (2005), no.~2, 239--252. 

\bibitem{3}C.~Petronio, A.~Zanellati: Algorithmic simplification of knot diagrams: new moves and experiments,
J. Knot Theory Ramifications 25 (2016), no.~10, 1650059, 30~pp. 

\bibitem{4}K. Reidemeister, "Elementare Begrundung der Knotentheorie" Abh. Math. Sem. Univ. Hamburg, 5 (1927), pp.24--32.

\bibitem{5}M.~Khovanov:
Doodle groups,
Trans. Amer. Math. Soc. 349 (1997), no.~6, 2297--2315. 

\bibitem{6}Y.~Miyazawa: A distance for diagrams of a knot, Topology Appl. 159 (2012), no.~4, 1122--1131.

\bibitem{7}Z.~Cheng, H.~Gao: A note on the independence of Reidemeister moves,
J. Knot Theory Ramifications 21 (2012), no.~9, 1220001, 7 pp. 

                               
\end{thebibliography}
\end{document}